\documentclass[10pt,leqno]{article}

\usepackage{amsthm,amsmath,amssymb}
\usepackage{dsfont}
\newcommand{\bfm}[1]{\boldsymbol{#1}}

\def \H{ \mathbb H}
\def \M{\mathbb M}
\def \R{ \mathbb R }
\def \reel{ \mathbb R }
\def \nat{ \mathbb N}

\def \ent{\mathds Z}
\def \one{ {\rm 1}\mkern-4.5mu{\rm I} }

\def \E{ \mathbb  E }

\def \P{ \mathbb P  }

\newtheorem{theorem}{Theorem}[section]

\newtheorem{definition}[theorem]{Definition}
\newtheorem{lemma}[theorem]{Lemma}
\newtheorem{proposition}[theorem]{Proposition}

\newtheorem*{remark*}{Remark}

\newtheorem*{remarkopen}{Remark and open problem}
\begin{document}

% "Title of the paper"
\title{Stochastic Homogenization of
Reflected Stochastic Differential Equations}
%\runtitle{Stochastic Homogenization of RSDE}

% indicate corresponding author with \corref{}
 \author{R\'emi Rhodes\\Ceremade-Universit\'e Paris-Dauphine\\ Place du mar\'echal De Lattre de Tassigny\\ 75775 Paris cedex 16\\
 email: \textsf{rhodes@ceremade.dauphine.fr}}
% \thankstext{t1}{Thanks} 
\date{\empty}

%\author{\fnms{cc} \snm{???}\ead[label=e1]{???}}
%\address{\printead{e1}}
%\and
%\author{\fnms{???} \snm{???}\ead[label=e2]{???}}
%\address{\printead{e2}}
%\affiliation{???}
\maketitle

\begin{abstract}
We investigate a functional limit theorem (homogenization) for Reflected Stochastic
Differential Equations on a half-plane with stationary coefficients when it is necessary to analyze both the effective Brownian motion and the effective local time. We prove that the limiting process is a reflected non-standard Brownian motion. Beyond the result, this problem is known as a prototype of non-translation invariant problem making the usual  method of the "environment as seen from the particle" inefficient.
\end{abstract}

\begin{center}
AMS: 60K37, 60F17,74Q99\\
Keywords: functional limit theoem, reflected Brownian motion, random medium, Skorohod problem, local time.
\end{center}

%%%%%%%%%%%%%%%%%%%%%

\section{Introduction}\label{sec:intro}
%%%%%%%%%%%%%%%%%%%%%%%%%%%%%%%%%%%%%%%%%%%%%%%%%%%%%%%%%%%%%%%%%%%%%%%%%%%%%%%%%%%%%%%%%%%%%%%%%%%%%%%%%%%%%%%%%%%%%%%%%%%
\subsection*{Statement of the problem} This paper is concerned with homogenization of Reflected
Stochastic Differential Equations (RSDE for short) evolving in a random medium, that is (see e.g.
\cite{jikov})
\begin{definition}{(Random medium).}\label{medium}
Let $(\Omega ,{\cal G},\mu )$ be a probability space and $\left\{\tau_{x};x\in \reel ^d\right\}$ be a group of measure
preserving transformations acting ergodically on $\Omega $, that is:

1) $\forall A\in {\cal G},\forall x\in \reel^d$, $\mu (\tau
_{x}A)=\mu (A)$,

2) If for any $x\in \reel^d$, $\tau _{x}A=A$ then $\mu (A)=0$ or
$1$,

3) For any measurable function ${\bfm g}$ on $(\Omega ,{\cal
G},\mu)$, the function $(x,\omega )\mapsto {\bfm g}(\tau_x \omega)$
is measurable on $(\reel^d\times\Omega ,{\cal B}(\reel^d)\otimes
{\cal G})$.
\end{definition}
The expectation with respect to the random medium is denoted by
${\mathbb M}$. In what follows we shall use the bold type to denote a
random function ${\bfm g}$ from $\Omega\times \R^p$ into $\R^n$ ($n \geq 1$ and $p\geq 0$).

A random medium is a mathematical tool to define stationary random functions. Indeed, given a function ${\bfm f}:\Omega \to \R$, we can consider for each fixed $\omega$ the function $x\in\R^d\mapsto {\bfm f}(\tau_x\omega)$.  This is a random function (the parameter $\omega$ stands for the randomness) and because of 1) of Definition \ref{medium}, the law of that function is invariant under $\R^d$-translations, that is both functions $ {\bfm f}(\tau_\cdot \omega)$ and $ {\bfm f}(\tau_{y+\cdot} \omega)$ have the same law for any $y\in\R^d$. For that reason, the random function is said to be stationary.

We suppose that we are given a random $d\times d$-matrix valued function ${\bfm \sigma}:\Omega\to \R^{d\times d} $, two random vector valued functions ${\bfm b},{\bfm \gamma}:\Omega\to \R^{d}$ and a
d-dimensional Brownian motion $B$ defined on a complete probability space  $(\Omega',{\cal F},\P)$ (the Brownian motion and the random medium
are independant). We shall
describe the limit in law, as $\varepsilon$ goes to $0$, of the following RSDE with stationary random coefficients 
\begin{equation}\label{SDEintro}
\begin{split}
dX^\varepsilon_t&=\varepsilon^{-1}{\bfm b}(\tau_{X^\varepsilon_t/\varepsilon}\omega)\,dt
+{\bfm \sigma}(\tau_{X^\varepsilon_t/\varepsilon}\omega)\,dB_t+{\bfm \gamma}(\tau_{X^\varepsilon_t/\varepsilon}\omega)
\,dK^\varepsilon_t,
\end{split}
\end{equation}
where $X^\varepsilon,K^\varepsilon$ are $({\cal F}_t)_t$-adapted processes ($ {\cal F}_t$ is the $\sigma$-field generated by $B$ up to time $t$) with constraint $X^\varepsilon_t\in\bar{D}$, where $D\subset\R^d$ is
the half-plane $\{(x_1,\dots,x_d)\in \R^d; x_1>0\}$,  $K^\varepsilon$ is the so-called
local time of the process $X^\varepsilon$, namely a continuous
nondecreasing process, which only increases on the set
$\{t;X^\varepsilon_t\in \partial D\}$. The reader is referred to  \cite{lions}  for strong existence and uniqueness results to \eqref{SDEintro} (see e.g \cite{watanabe} for the weak existence), in particular under the assumptions on the coefficients ${\bfm \sigma},{\bfm b}$ and  ${\bfm \gamma}$  listed below. Those stochastic processes
are involved in the probabilistic representation of second order
 partial differential equations in half-space with Neumann
boundary conditions (see \cite{pardouxzhang} for an insight of the topic). In particular, we are interested in
homogenization problems for which it is necessary to identify both
the homogenized equation and the homogenized boundary conditions.

Without the reflection term ${\bfm \gamma}(X^\varepsilon_t/\varepsilon)
\,dK^\varepsilon_t$, the issue of determining the limit in
\eqref{SDEintro} is the subject of an extensive literature in the
case when the coefficients ${\bfm b},{\bfm \sigma}$ are periodic, quasi-periodic
and, more recently, evolving in a stationary ergodic random medium.
Quoting all references is beyond the scope of this paper. Concerning
 homogenization of RSDEs, there are only a few works dealing with
periodic coefficients (see \cite{arisawa,barles,bensoussan,tanaka}).
As pointed out in \cite{barles}, homogenizing \eqref{SDEintro} in a random medium is a well-known problem that remains unsolved yet. There are several
difficulties in this framework that make the classical machinery of
diffusions in random media (i.e. without reflection) fall short of
determining the limit in \eqref{SDEintro}. In particular, the reflection term
breaks the stationarity properties of the process $X^\varepsilon$ so
that the method of the {\it environment as seen from the particle}
(see \cite{olla} for an insight of the topic) is inefficient. Moreover, the lack of compactness of a random medium prevents from using compactness methods. The main resulting difficulties are the lack of invariant probability measure (IPM for short) associated to the process $X^\varepsilon$ and the study of the boundary ergodic problems. The aim of this paper is precisely to investigate the random case and prove the convergence of the process $X^\varepsilon$ towards a reflected Brownian motion. The convergence is established in probability with respect to the random medium and the starting point $x$. 

We should also point out that the problem of determining the limit in \eqref{SDEintro} could be expressed in terms of reflected random walks in random environment, and remains quite open as well. In that case, the problem could be stated as follows: suppose we are given, for each $z\in \ent^d$ satisfying $|z|=1$, a random variable ${\bfm c}(\cdot,z):\Omega\to ]0;+\infty[$. Define the continuous time process $X$ with values in the half-lattice $L=\nat\times \ent^{d-1}$ as the random walk that, when arriving at a site $x\in L$, waits a random exponential time of parameter $1$ and then performs a jump to the neighboring sites $y\in L$ with jump rate ${\bfm c}(\tau_x\omega,y-x)$. Does the rescaled random walk $\varepsilon X_{t/\varepsilon^2}$ converge in law towards a reflected Brownian motion? Though we don't treat explicitly that case, our proofs can be adapted to that framework.

\subsection*{Structure of the coefficients}

\textit{Notations:} Throughout the paper,
we use the convention of summation over repeated indices
$\sum_{i=1}^d c_id_i=c_id_i $ and we use the superscript $\phantom{}^*$ to denote the transpose $A^*$ of some given matrix $A$. If a random function ${\bfm \varphi}:\Omega\to \R$ possesses smooth trajectories, i.e.  for any $\omega\in \Omega$ the mapping  $x\in\R^d\mapsto {\bfm \varphi}(\tau_x\omega) $ is smooth with bounded derivatives, we can consider
its partial derivatives at $0$ denoted by $D_i{\bfm \varphi}$, that is $D_i{\bfm \varphi}(\omega)=\partial_{x_i}(x\mapsto {\bfm \varphi}(\tau_x\omega) )_{|x=0}$. \qed

\vspace{2mm}
We define ${\bfm a}={\bfm \sigma}{\bfm \sigma}^*$. %We assume that we are given an antisymmetric matrix-valued
%function ${\bfm H}:\Omega\rightarrow
%\R^{d\times d}$ such that  ${\bfm H}_{ij}=0$ whenever $i=1$ or $j=1$. 
For the sake of simplicity, we assume that $\forall \omega\in \Omega$ the mapping $x\in\R^d\mapsto {\bfm \sigma}(\tau_x\omega) $ %, $x\mapsto {\bfm H}(\tau_x\omega)$ are
is bounded and smooth with bounded derivatives of all orders. We further impose these bounds do not depend on $\omega$. 

Now we motivate the structure we impose on the coefficients ${\bfm b}$ and ${\bfm \gamma}$. A specific point in the literature of diffusions in random media is that the lack of compactness of a random medium makes it impossible to find an IPM for the involved diffusion process.
There is a simple argument to understand why: since the coefficients of the SDE driving the $\R^d$-valued diffusion process are stationary, any $\R^d$-supported invariant measure must be stationary. So, unless it is trivial, it cannot have finite mass. That difficulty has been overcome by introducing the "environment as seen from the particle" (ESFP for short). It is a $\Omega$-valued Markov process describing the configurations of the environment visited by the diffusion process: briefly, if you denote by $X$ the diffusion process then the ESFP should match $\tau_{X}\omega$. There is a well known formal ansatz that says: if we can find a bounded function ${\bfm f}:\Omega\to [0,+\infty[$ such that, for each $\omega\in \Omega$, the measure ${\bfm f}(\tau_{x }\omega)dx $ is invariant for the diffusion process, then the probability measure ${\bfm f}(\omega)d\mu$ (up to a renormalization constant) is invariant for the ESFP. So we can switch an invariant measure with infinite mass asociated to the diffusion process for an IPM associated to the ESFP. 

The remaining problem is  to find an invariant measure (of the type ${\bfm f}(\tau_x\omega)dx $) for the diffusion process. Generally speaking, there is no way to find it  excepted when it is explicitly known. In the stationary case (without reflection), the most general situation when it is explicitly known is when the generator of the rescaled diffusion process can be rewritten in divergence form as 
\begin{equation}\label{lineargen}
{\cal L}^\varepsilon f=\frac{1}{2}e^{2{\bfm V}(\tau_{x/\varepsilon}\omega)}\partial_{x_i}\big(e^{-2{\bfm V}(\tau_{x/\varepsilon}\omega)}({\bfm a}_{ij}+{\bfm H}_{ij})(\tau_{x/\varepsilon}\omega)\partial_{x_j}f\big),
\end{equation}
where ${\bfm V}:\Omega\to\R$ is a bounded scalar function and ${\bfm H}:\Omega\to \R^{d\times d}$ is a function taking values in the set of antisymmetric matrices. The invariant measure is then given by $e^{2{\bfm V}(\tau_{x/\epsilon}\omega)}dx $ and the IPM for the ESFP matches $e^{2{\bfm V}(\omega)}d\mu$. However, it is common  to assume ${\bfm V}={\bfm H}=0$ to simplify the problem since the general case is in essence very close to that situation. Why is the existence of an IPM so important? Because it entirely  determines the asymptotic behaviour of the diffusion process via ergodic theorems. The ESFP is therefore a central point in the literature of diffusions in random media.

The case of RSDE in random media does not derogate this rule and we are bound to find a framework where the invariant measure is (at least formally) explicitly known. So we assume that the entries of the coefficients ${\bfm b}$ and ${\bfm \gamma} $, defined on $\Omega$, are given by
\begin{align}\label{defcoef}
 \forall j =1,\dots ,d,\,\,\,{\bfm b}_j=\frac{1}{2}D_i{\bfm a}_{ij},\quad {\bfm \gamma}_j =  {\bfm a}_{j1}.%+D_i {\bfm H}_{ji}.
\end{align}
With this definition, the generator of the Markov process $X^\varepsilon$ can be rewritten in divergence form as (for a sufficiently smooth function $f$ on $\bar{D}$) 
\begin{equation}\label{gen:intro}
{\cal
L}^\varepsilon f=\frac{1}{2}\partial_{x_i}\big({\bfm a}_{ij}(\tau_{x/\varepsilon}\omega)\partial_{x_j}f\big)
\end{equation}
with boundary condition ${\bfm \gamma}_i(\tau_{x/\varepsilon}\omega)\partial_{x_i}f=0$ on $\partial D$. If the environment $\omega$ is fixed, it is a simple exercise to check that the Lebesgue measure is formally invariant for the process $X^\epsilon$. If the ESFP exists, the aforementioned ansatz tells us that  $\mu$ should be an IPM for the ESFP. Unfortunately, we shall see that there is no way of defining properly the ESFP. The previous formal discussion is however helpful to provide a good intuition of the situation and to figure out what the correct framework must be. Furthermore the framework \eqref{defcoef} also comes from physical motivations. As defined above, the reflection term ${\bfm \gamma}$ coincides with the so-called conormal field and the associated PDE problem is said to be of Neumann type. From the physical point of view, the conormal field is the "canonical" assumption that makes valid the mass conservation law since the relation ${\bfm a}_{j1}(\tau_{x/\varepsilon}\omega)\partial_{x_j}f=0$ on $ \partial D$ means that the flux through the boundary must vanish.
 Our framework for RSDE is therefore to be seen as a natural generalization of the classical stationary framework. 

\begin{remark*}
It is straightforward to adapt our proofs to treat the more general situation when the generator of the RSDE inside $D$ coincides with \eqref{lineargen}. In that case, the reflection term is given by $ {\bfm \gamma}_j =  {\bfm a}_{j1}+{\bfm H}_{j1}  $. 
\end{remark*}

% It is then natural to wonder if we can consider a larger class of reflection coefficients ${\bfm \gamma}$ while preserving the mass conservation law. The answer is positive if ${\bfm \gamma}$ has the structure described in \eqref{defcoef} and this remark motivates our definition of ${\bfm \gamma}$.

Without loss of generality, we assume that $a_{11}=1$. We  further assume that  ${\bfm a}$ is uniformly elliptic, i.e. there exists a constant $\Lambda>0$ such that
\begin{equation}\label{unif_ellip}
\forall\omega\in \Omega,\quad \Lambda{\rm I}\leq {\bfm
a}(\omega)\leq \Lambda^{-1} {\rm I}.
\end{equation}
That assumption means that the process $X^\epsilon$ diffuses enough, at each point of $\bar{D}$, in all directions. It is thus is a convenient assumption to ensure the ergodic properties of the model. The reader is referred, for instance, to \cite{delarue,rhodes:APS,sznitman} for various situations going beyond that assumption. We also point out that, in the context of RSDE, the problem of homogenizing \eqref{SDEintro} without assuming \eqref{unif_ellip} becomes quite challenging, especially when dealing with the boundary phenomena.

\subsection*{Main Result} In what follows, we indicate by $\P_x^\varepsilon$ the law of
the process $X^\varepsilon$ starting from $x\in\bar{D}$ (keep in mind that this probability measure also depends on $\omega$ though it does not appear through the notations). Let us consider  a  nonnegative function $\chi:\bar{D}\to \R_+$ such that $\int_{\bar{D}}\chi(x)\,dx=1$. Such a function defines a probability measure on $\bar{D}$ denoted by $\chi(dx)=\chi(x)dx$. We fix $T>0$. Let $C$ denote the space of continuous $\bar{D}\times\R_+ $-valued functions on $[0,T]$ equipped with the sup-norm topology. We are now in position to state the main result of the paper:
\begin{theorem}\label{mainth}
The C-valued process $(X^\varepsilon,K^\varepsilon)_\varepsilon$ converges weakly, in
$\mu\otimes \chi$ probability, towards the solution
$(\bar{X},\bar{K})$ of the RSDE
\begin{equation}\label{SDElimit}
\bar{X}_t=x+\bar{A}^{1/2}B_t+\bar{\Gamma}\bar{K}_t,
\end{equation} with constraints
$\bar{X}_t\in\bar{D}$ and $\bar{K}$ is the local time associated to $\bar{X}$. In other words, for each bounded continuous function $F$ on $C$ and $ \delta>0$, we have
\begin{equation*}
\lim_{\varepsilon\to 0}\,\,\,\mu\otimes \chi \left\{(\omega,x)\in \Omega\times \bar{D};\big|\E^\varepsilon_x(F(X^\varepsilon,K^\varepsilon))-\E_x(F(\bar{X},\bar{K}))\big|\geq \delta\right\}=0.
\end{equation*}
The so-called homogenized (or effective) coefficients $\bar{A}$ and $\bar{\Gamma}$ are constant. Moreover $\bar{A}$ is invertible, obeys a variational formula (see subsection \ref{corr} for the meaning of the various terms)
$$\bar{A}=\inf_{{\bfm \varphi}\in\mathcal{C}}\M\big[({\rm I}+D{\bfm \varphi})^*{\bfm a}({\rm I}+D{\bfm \varphi})\big],$$ and $\bar{\Gamma}$ is the conormal field associated to $\bar{A}$, that is $\bar{\Gamma}_i=\bar{A}_{1i}$ for $i=1,\dots,d$.
\end{theorem}

\begin{remarkopen}
The reader may wonder whether it may be simpler to consider the case ${\bfm \gamma}_{i}=\delta_{1i}$ where $\delta$ stands for the Kroenecker symbol. In that case, ${\bfm \gamma}$ coincides with the normal to $\partial D$. Actually, this situation is much more complicated since one can easily be convinced that there is no obvious invariant measure associated to $X^\epsilon$.

On the other side, one may wonder if, given the form of the generator \eqref{gen:intro} inside $D$, one can find a larger class of reflection coefficients ${\bfm \gamma}$ for which the homogenization procedure can be carried through. Actually, a computation based on the Green formula shows that it is possible to consider a bounded antisymmetric matrix valued function ${\bfm A}:\Omega\to \R^{d\times d}$ such that  ${\bfm A}_{ij}=0$ whenever $i=1$ or $j=1$, and to set $ {\bfm \gamma}_j =  {\bfm a}_{j1}+D_i {\bfm A}_{ji}$. In that case, the Lebesgue measure is invariant for $X^\epsilon $. Furthermore, the associated Dirichlet form (see subsection \ref{ergo}) satisfies a strong sector condition in such a way that the construction of the correctors is possible. However, it is not clear whether the localization technique based on the Girsanov transform (see Section \ref{secgirsanov} below) works. So we leave that situation as an open problem.
\end{remarkopen}

The non-stationarity of the problem makes the proofs technical. So we have  divided the remaining part of the paper  into  two parts. In order to have a global understanding  of the proof of Theorem \ref{mainth}, we set the main steps out  in Section \ref{guide}  and gather most of the technical proofs in the Appendix.

%%%%%%%%%%%%%%%%%%%%%%%%%%%%%%%%%%%%%%%%%%%%%%%%%%%%%%%%%%%%%%%%%%%%%%%%%%%%%%%%%%%%%%%%%%%%%%%%%%%%%%%%%%%%%%%%%%%%%%%%%%%

\section{Guideline of the proof}\label{guide}

%%%%%%%%%%%%%%%%%%%%%%%%%%%%%%%%%%%%%%%%%%%%%%%%%%%%%%%%%%%%%%%%%%%%%%%%%%%%%%%%%%%%%%%%%%%%%%%%%%%%%%
As explained in introduction, what makes the problem of homogenizing RSDE in random medium known as a difficult problem is the lack of stationarity of the model. The first resulting difficulty is that you cannot define properly the ESFP (or a somewhat similar process) because you cannot prove that it is a Markov process. Such a process is important since its IPM encodes what the asymptotic behaviour of the process should be. The reason why the ESFP is not a Markov process is the following. Roughly speaking, it stands for an observer sitting on a particle $X^\epsilon_t$ and looking at the environment $\tau_{X^\epsilon_t}\omega $ around the particle. For this process to be Markovian, the observer must be able to determine, at a given time $t$, the future evolution of the particle with the only knowledge of the  environment $\tau_{X^\epsilon_t}\omega $. In the case of RSDE, whenever the observer sitting on the particle wants to determine the future evolution of the particle, the knowledge of the environment $\tau_{X^\epsilon_t}\omega $ is not sufficient. He  also needs to know whether the particle is located on the boundary $\partial D$ to determine if the pushing of the local time $K^\epsilon_t$ will affect the trajectory of the particle. So we are left with the problem of dealing with a process $X^\epsilon$ possessing no IPM.

%Second, the non-stationarity of the model will affect the construction of the so-called correctors involved in the ergodic problems for the local time. It is well known that, to be efficient, the correctors must satisfy a sublinear growth condition, which is usually derived from the stationarity properties of the model. In our context, those used to determine the limit as $\epsilon\to 0 $ of the functional 
%$$\int_0^t{\bfm f}(\tau_{X^\epsilon_r/\epsilon}\omega)\,dK^\epsilon_r$$
%for some generic function ${\bfm f}:\Omega\to \R$ exhibit a bad behaviour at infinity due to the lack of stationarity of the problem.

\subsection{Localization}\label{secgirsanov}

To overcome the above difficulty, we shall use a localization technique. Since the process $X^\epsilon$ is not convenient to work with, the main idea is to compare  $X^\epsilon$ with a better process that possesses, at least locally, a similar asymptotic behaviour. To be better, it must have an explicitly known IPM. There is a simple way to find such a process: we plug a smooth and deterministic potential $V:\bar{D}\to \R$  into \eqref{gen:intro} and define a new operator acting on $ C^2(\bar{D})$
\begin{equation}\label{gen:introV}
{\cal
L}^\varepsilon_V=\frac{e^{2V(x)}}{2}\sum_{i,j=1}^d\partial_{x_i}\big(e^{-2V(x)}{\bfm a}_{ij}(\tau_{x/\varepsilon}\omega)\partial_{x_j}\,\big)={\cal
L}^\varepsilon-\partial_{x_i}V(x){\bfm a}_{ij}(\tau_{x/\varepsilon}\omega)\partial_{x_j},
\end{equation}
with the same boundary condition ${\bfm \gamma}_{i}(\tau_{x/\varepsilon}\omega)\partial_{x_i}=0 $ on $\partial D $. If we impose the condition 
\begin{equation}\label{V1}
 \int_{\bar{D}}e^{2V(x)}dx=1
 \end{equation}
  and fix the environment $\omega$, we shall prove that the RSDE with generator ${\cal L}^\varepsilon_V$ inside $D$ and boundary condition ${\bfm \gamma}_{i}(\tau_{x/\varepsilon}\omega)\partial_{x_i}=0 $ on $\partial D $ admits $e^{2V(x)}dx$ as IPM. 
  
  Then we want to find a connection between the process $X^\epsilon$ and the Markov process with generator ${\cal L}^\varepsilon_V$ inside $D$ and boundary condition ${\bfm \gamma}_{i}(\tau_{x/\varepsilon}\omega)\partial_{x_i}=0 $ on $\partial D $. To that purpose, we use the Girsanov transform. 
 % To compare the law of the process $X^\epsilon$ with that of the RSDE associated to ${\cal L}^%\varepsilon_V$, 
More precisely, we fix $T>0$ and impose
\begin{equation}\label{V2}
V \text{ is smooth and }\partial_xV \text{ is bounded}.
\end{equation} 
Then we define the following probability measure on the filtered space $(\Omega';\mathcal{F},(\mathcal{F}_t)_{0\leq t\leq T}) $ 
$$d\P^{\varepsilon*}_x=\exp\Big(-\int_0^T\partial_{x_i}V(X^\varepsilon_r){\bfm \sigma}_{ij}(\tau_{X^\varepsilon_r/\varepsilon}\omega)\,dB^j_r-\frac{1}{2}\int_0^T\partial_{x_i}V(X^\varepsilon_r){\bfm a}_{ij}(\tau_{X^\varepsilon_r/\varepsilon}\omega)\partial_{x_j}V(X^\varepsilon_r)\,dr\Big)\,d\P^\varepsilon_x.$$
Under $\P^{\varepsilon*}_x$,  the process $B^*_t=B_t+\int_0^t{\bfm \sigma}(\tau_{X^\varepsilon_r/\varepsilon}\omega)\partial_x V(X^\varepsilon_r)\,dr$ ($0\leq t\leq T$) is a Brownian motion and  the process $X^\varepsilon$ solves the RSDE
\begin{equation}\label{SDEM}
dX^\varepsilon_t=\varepsilon^{-1}{\bfm b}(\tau_{X^\varepsilon_t/\varepsilon}\omega)\,dt-{\bfm a}(\tau_{X^\varepsilon_t/\varepsilon}\omega)\partial_x V(X^\varepsilon_t)\,dt+{\bfm \sigma}(\tau_{X^\varepsilon_t/\varepsilon}
\omega)\,dB_t^*+{\bfm \gamma}(\tau_{X^\varepsilon_t/\varepsilon}\omega)
\,dK^\varepsilon_t
\end{equation}
starting from $X^\varepsilon_0=x$, where $K^\varepsilon$ is the local time of $X^\varepsilon$.  It is straightforward to check that, if $B^*$ is a Brownian motion, the generator associated to the above RSDE coincides with \eqref{gen:introV} for sufficiently smooth functions. To sum up, with the help of the Girsanov transform, we can compare the law of the process $X^\epsilon$ with that of the RSDE \eqref{SDEM} associated to ${\cal L}^\varepsilon_V$.

We shall see that most of the necessary estimates to homogenize  the process $X^\epsilon$ are valid under $\P^{\varepsilon*}_x$. We want to make sure that they remain valid under $\P^{\varepsilon}_x$. To that purpose, the probability measure $\P^{\varepsilon}_x$ must be dominated by $\P^{\varepsilon*}_x$ uniformly with respect to $\epsilon$. From \eqref{V2}, it is readily seen that $C=\sup_{\epsilon>0}\big(\E^{\varepsilon*}_x\big[(\frac{d\P^\epsilon_x}{d\P^{\epsilon*}_x})^2\big]\big)^{1/2}<+\infty$ ($C$ only depends on $T,|{\bfm a}|_{\infty}$ and $\sup_{\bar{D}}|\partial_xV|$). Then the  Cauchy-Schwarz inequality yields
\begin{equation}\label{girsanov}
\forall\epsilon>0,\,\,\forall A\,\,\, \mathcal{F}_T\text{-measurable subset },\quad \P^\varepsilon_x(A)\leq C\big(\P^{\varepsilon*}_x(A)\big)^{1/2}.
\end{equation}

In conclusion, we summarize our strategy: first we shall prove that the process $X^\epsilon$ possesses an IPM under the modified law $\P^{\varepsilon*}$, then we establish under $\P^{\varepsilon*}$ all the necessary estimates to homogenize $X^\epsilon$ , and finally we shall deduce that the estimates remain valid under $\P^\varepsilon$ thanks to \eqref{girsanov}. Once that is done, we shall be in position to homogenize \eqref{SDEintro}.

To fix the ideas and to see that the class of functions $V$ satisfying \eqref{V1} \eqref{V2}  is not empty, we can choose $V$ to be equal to 
\begin{equation}\label{defV}
V(x_1,\dots,x_d)=Ax_1+A(1+x_2^2+\dots+x_d^2)^{1/2}+c,
\end{equation} 
for some renormalization constant $c$ such that $\int_{\bar{D}}e^{-2V(x)}\,dx=1$ and some positive constant $A$.

%%%%%%%%%%%%%%%%%%%%%%%%%
%%%%%%%%%%%%%%%%%%%%%%%%%%%%%%%%

\vspace{2mm}

 \textit{Notations for measures.}  In what follows, $\bar{\P}^\varepsilon$ (resp. $\bar{\P}^{\varepsilon*}$) stands for
the  averaged (or annealed) probability measure
$\M\int_{\bar{D}}\P_x^{\varepsilon}(\cdot)e^{-2V(x)}\,dx$ (resp. $\M\int_{\bar{D}}\P_x^{\varepsilon*}(\cdot)e^{-2V(x)}\,dx$), and
$\bar{\E}^\varepsilon$ (resp. $\bar{\E}^{\varepsilon*}$) for the corresponding expectation.  \\
$\P_D^*$ and $\P_{\partial D}^* $ respectively denote the probability measure $e^{-2V(x)}\,dx\otimes d\mu$ on $\bar{D}\times \Omega$ and the finite measure $e^{-2V(x)}\,dx\otimes d\mu$ on $\partial D \times \Omega$. $\M_D^*$ and $\M_{\partial D}^*$ stand for the respective expectations.

\subsection{Invariant probability measure}\label{ipm}
%%%%%%%%%%%%%%%%%%%%%%%%%%%%%%%%%%%%%%%%%%%
As explained above, the main advantage of considering the process $X^\epsilon$ under the modified law $\P^{\varepsilon*}_x$ is that we can find an IPM. More precisely
\begin{lemma}\label{lem_invariant}The process $X^\varepsilon$ satisfies:\\
1)For each function ${\bfm f}\in L^1(\bar{D}\times \Omega
;\P_D^*)$ and $t\geq 0$:
\begin{equation}\label{invariant}
\bar{\E}^{\varepsilon *}[{\bfm
f}(X^{\varepsilon}_t,\tau_{X^{\varepsilon}_t/\varepsilon}\omega
)]=\M^* _{D}[{\bfm f}].
\end{equation}
2) For each function ${\bfm f}\in L^1(\partial D \times \Omega
;\P^*_{\partial D})$  and $t\geq 0$:
\begin{equation}\label{invariant_local}
\bar{\E}^{\varepsilon *} \big[\int_0^t{\bfm
f}(X^\varepsilon_r,\tau_{X^\varepsilon_r/\varepsilon}\omega)\,dK^\varepsilon_r\big]
= t\M_{\partial D} ^*\big[ {\bfm f}\big].
\end{equation}
\end{lemma}

The first relation \eqref{invariant} results from the structure of $\mathcal{L}^\epsilon_V$ (see \eqref{gen:introV}), which has been defined so as to make $e^{-2V(x)}dx$ invariant for the process $X^\epsilon$. Once \eqref{invariant} established, \eqref{invariant_local} is derived from the fact that $K^\epsilon$ is the density of occupation time of the process $X^\epsilon$ at the boundary $\partial D$.

\subsection{Ergodic problems}\label{ergo}
%%%%%%%%%%%%%%%%%%%%%%%%%%%%%%%%%%%%%%%%%%%

The next step is to determine the asymptotic behaviour as $\epsilon\to 0 $ of the quantities
\begin{equation}\label{ergfunc}
\int_0^t{\bfm f}(\tau_{X^\epsilon_r/\epsilon}\omega)\,dr\quad \text{ and }\quad \int_0^t{\bfm f}(\tau_{X^\epsilon_r/\epsilon}\omega)\,dK^\epsilon_r. 
\end{equation}
The behaviour of each above quantity is related to the evolution of the process $X^\epsilon$ respectively inside the domain $D$ and near the boundary $\partial D$. We shall see that both limits  can be identified  by solving ergodic problems associated to appropriate resolvent families. What concerns the first functional has already been investigated in the literature. The main novelty of the following section is the boundary ergodic problems associated to the second functional.

\subsubsection*{Ergodic problems associated to the diffusion process inside $D$}

First we have to figure out what happens when the process $X^\epsilon$ evolves inside the domain $D$. In that case, the pushing of the local time in \eqref{SDEintro} vanishes. The process $X^\epsilon$ is thus driven by the same effects  as in the stationary case (without reflection term). The ergodic properties of the process inside $D$ are therefore the same as in the classical situation. So we just sum up the main results and give references for further details.

\textit{Notations:}  For $p\in [1;\infty]$,  $L^p(\Omega)$ denotes the standard space of $p$-th power integrable functions (essentially bounded functions if $p=\infty$) on $(\Omega ,{\cal G},\mu )$ and $|\, \cdot \,|_p$ the corresponding norm. If $p=2$, the associated inner product is denoted by $(\,
\cdot\,,\, \cdot\,)_2$.  The space $C^\infty_c(\bar{D})$ (resp. $C^\infty_c(D)$) denotes the space of smooth functions on $\bar{D}$ with compact support in $\bar{D}$ (resp. $D$). \qed

\textit{Standard background:} The operators on $L^2(\Omega )$ defined by
$T_{x}{\bfm g}(\omega )={\bfm g}(\tau_{x}\omega )$ form a strongly
continuous group of unitary maps in $L^2(\Omega )$. Let $(e_1,\dots,e_d)$ stand
for the canonical basis of $\R^d$. The group $(T_x)_x$ possesses $d$
generators defined by $D_i{\bfm g}=\lim_{ h\in\R\to
0}h^{-1}(T_{he_i}{\bfm g}-{\bfm g})$, for  $i=1,\dots,d$,   whenever the limit exists in
the $L^2(\Omega)$-sense. The operators $(D_i)_i$ are closed and densely defined. Given ${\bfm \varphi}\in \bigcap_{i=1}^d{\rm Dom}(D_i)$, $D{\bfm \varphi}$ stands for the d-dimensional vector whose entries are $D_i{\bfm \varphi}$ for $i=1,\dots,d$. 

We point out that we distinguish $D_i$ from the usual differential operator $\partial_{x_i}$ acting on differentiable functions $f:\R^d\to \R$ (more generally, for $k\geq 2$, $\partial^k_{x_{i_1}\dots x_{i_k}}$ denotes the iterated operator $ \partial_{x_{i_1}}\dots \partial_{x_{ik}}$). However, it is straightforward to check that,  whenever a function ${\bfm \varphi}\in {\rm Dom}(D)$ possesses differentiable trajectories (i.e. $\mu$ a.s. the mapping $x \mapsto {\bfm \varphi}(\tau_x\omega)$ is differentiable in the classical sense), we have $D_i{\bfm \varphi}(\tau_x\omega)=\partial_{x_i} {\bfm \varphi}(\tau_x\omega) $. \\

We denote by ${\cal C}$ the dense subspace of $L^2(\Omega )$ defined by
\begin{equation}\label{defc}
{\cal C}={\rm Span}\left\{{\bfm g} \star \varphi ;{\bfm g}\in
L^\infty(\Omega ),\varphi \in C^\infty _c(\reel ^{d} )\right\}\quad \text{ where }   {\bfm g} \star \varphi(\omega )=\int_{\reel ^{d}}{\bfm
g}(\tau_{x}\omega )\varphi (x)\,dx
\end{equation}
Basically, ${\cal C}$ stands for the space of smooth functions on the random medium. We have ${\cal
C}\subset {\rm Dom}(D_i)$  and $D_i({\bfm
g} \star \varphi)=-{\bfm g} \star
\partial_{x_i}  \varphi$ for all $1 \leq i \leq d$. This quantity is also equal to $D_i{\bfm g}
\star \varphi $ if ${\bfm g}\in {\rm Dom}(D_i)$. \qed

\vspace{2mm}

We associate to the operator ${\cal L}^{\varepsilon} $ (Eq.
\eqref{gen:intro}) an unbounded operator acting on ${\cal C}~\subset~
L^2(\Omega)$
\begin{equation}\label{gen_medium}
{\bfm L}=\frac{1}{2}D_i\big({\bfm a}_{ij}D_j\cdot \big).
\end{equation}
Following \cite[Ch. 3, Sect 3.]{fukushima} (see also \cite[Sect.
4]{rhodes:06}), we can consider its Friedrich extension, still
denoted by ${\bfm L}$, which is a self-adjoint operator on
$L^2(\Omega)$. The domain $\H$ of the corresponding Dirichlet form
can be described as the closure of ${\cal C}$ with respect to the
norm $\|{\bfm \varphi} \|_\H^2=|{\bfm \varphi}|_2^2+|D{\bfm
\varphi}|_2^2$. Since ${\bfm L}$ is self-adjoint, it also defines a
resolvent family $(U_\lambda)_{\lambda>0}$. For each ${\bfm f}\in
L^2(\Omega)$, the function ${\bfm w}_\lambda=U_\lambda({\bfm f})\in
\H\cap {\rm Dom}({\bfm L})$ equivalently solves the
$L^2(\Omega)$-sense equation
\begin{equation}\label{l2e}
\lambda{\bfm w}_\lambda-{\bfm L}{\bfm w}_\lambda={\bfm f}
\end{equation}
 or the weak formulation equation
\begin{equation}\label{wfe}
  \forall {\bfm \varphi}\in \H,\quad \lambda({\bfm w}_\lambda,{\bfm
\varphi})_2+(1/2)\big({\bfm a}_{ij}D_i{\bfm w}_\lambda,D_j{\bfm
\varphi}\big)_2=({\bfm f},{\bfm \varphi})_2.
\end{equation}
Moreover, the resolvent operator $U_\lambda$ satisifes the maximum principle: 
\begin{lemma}\label{maxU}
For any function ${\bfm f}\in L^\infty(\Omega)$, the function $U_\lambda({\bfm f})$ belongs to  $L^\infty(\Omega)$ and satisfies $$|U_\lambda({\bfm f})|_\infty\leq  |{\bfm f}|_\infty/\lambda.$$
\end{lemma}
The ergodic properties of the operator ${\bfm L}$ are summarized in the following proposition:
\begin{proposition}\label{ergprob}
Given ${\bfm f}\in L^2(\Omega)$, the solution ${\bfm w}_\lambda$ of
the resolvent equation  $  \lambda {\bfm
w}_\lambda-{\bfm L}{\bfm w}_\lambda={\bfm f}$ ($\lambda>0$) satisfies
\begin{equation*}
  |\lambda{\bfm w}_\lambda-\M[{\bfm f}]|_2\rightarrow 0\, \text{ as }\lambda \to
  0,\quad \text{and } \,\forall\lambda>0,\,\,|\lambda^{1/2}D{\bfm w}_\lambda|_2\leq \Lambda^{-1/2}|{\bfm
  f}|_2.
\end{equation*}
\end{proposition}

%%%%%%%%%%%%%%%%%%%%%%%%%%%%%%%%%%%%
\subsubsection*{Boundary ergodic problems}

Second, we have to figure out what happens when the process hits the boundary $\partial D$. If we want to adapt the arguments in \cite{tanaka}, it seems natural to look at the unbounded operator in random medium $H_\gamma $, whose construction is formally the following: given $ \omega\in\Omega$ and a smooth function
${\bfm \varphi}\in {\cal C}$, let us denote by
$\tilde{u}_\omega:\bar{D}\to \R$ the solution of the problem
\begin{equation}\label{uphi}
 \left\{\begin{array}{l}
L^\omega \tilde{u}_\omega(x)=0, \, x\in
  D,
  \\ \tilde{u}_\omega(x)={\bfm \varphi}(\tau_x\omega),\,x\in \partial D.
    \end{array}\right.
\end{equation}
where the operator $L^\omega$ is defined by
\begin{equation}\label{deflomega}
L^\omega
f(x)=(1/2){\bfm a}_{ij}(\tau_x\omega)\partial^2_{x_ix_j}f(x)+{\bfm b}_i(\tau_x\omega)\partial_{x_i}f(x)
\end{equation}
whenever $f:\bar{D}\to \R$ is smooth enough, say $f\in C^2(\bar{D})$. Then we define
\begin{equation}\label{defhgamma} H_\gamma{\bfm \varphi}(\omega)={\bfm
\gamma}_i(\omega)\partial_{x_i} \tilde{u}_\omega(0).
\end{equation}
\begin{remark*}Choose $\epsilon=1$ in \eqref{SDEintro} and denote by $(X^1,K^1)$ the solution of \eqref{SDEintro}. The operator $H_\gamma$ is actually the generator of the $\Omega $-valued Markov process $Z_t(\omega)=\tau_{Y_t(\omega)}\omega$, where $Y_t(\omega)=X^1_{K^{-1}(t)}$ and the function $K^{-1} $ stands for the left inverse of $K^1$: $K^{-1}(t)=\inf\{s>0;K^1_s\geq t\}$. The process $Z$ describes the environment as seen from the particle whenever  the process $X^1$  hits the boundary $\partial D$.
\end{remark*}

The main difficulty lies in constructing a unique solution of Problem
\eqref{uphi} with suitable growth and integrability properties
because of the lack of compactness of $D$. %
This point together with the lack of IPM are known as the major difficulties in homogenizing the Neumann problem in random media. We detail below the construction of $H_\gamma$ through
its resolvent family. In spite of its technical aspect, this contruction seems to be the right one because it exhibits a lack of stationarity  along the $e_1$-direction, which is intrinsec to the problem due to the pushing of the local time $K^\epsilon$, and conserves the stationarity of the problem along all other directions.

\vspace{2mm}
First we give a few notations before tackling the construction of $H_\gamma$. In what follows, the notation $(x_1,y)$ stands for a $d$-dimensional vector, where the first component $x_1$ belongs to $\R$ (eventually $\R_+=[0;+\infty)$) and the second component $y$ belongs to $\R^{d-1}$. To define an unbounded operator, we first need to determine the space that it acts on. As explained above, that space must exhibit a a lack of stationarity  along the $e_1$-direction and stationarity along all other directions. So the natural space to look at is the product space $\R_+\times \Omega$, denoted by $\Omega^+$, equipped with the measure $d\mu^+\stackrel{def}{=}dx_1\otimes d\mu$ where $dx_1$ is the Lebesgue measure on $\R_+$. We can then consider the standard spaces $L^p(\Omega^+)$ for $p\in[1;+\infty]$.

\vspace{2mm}
Our strategy is to define the Dirichlet form associated to $H_\gamma$. To that purpose, we need to define a dense space of test functions on $\Omega^+ $ and a symmetric bilinear form acting on the test functions. It is natural to define the space of test functions by
$$\mathds{C}(\Omega^+)={\rm Span}\{ \rho(x_1){\bfm \varphi}(\omega) ;\rho\in
C^\infty_c([0;+\infty)) ,{\bfm \varphi}\in \mathcal{C}\}.$$
Among the test functions we distinguish those that are vanishing on the boundary $\{0\}\times \Omega $ of $\Omega^+$
\begin{align*}
\mathds{C}_c(\Omega^+)&={\rm Span}\{ \rho(x_1){\bfm \varphi}(\omega) ;\rho\in
C^\infty_c((0;+\infty)) ,{\bfm \varphi}\in \mathcal{C}\}.
\end{align*}
 
Before tackling the construction of the symmetric bilinear form, we also need to introduce some elementary tools  of differential calculus on $\Omega^+$.
For any ${\bfm g}\in\mathds{C}(\Omega^+)$, we introduce a sort of gradient  $\partial {\bfm g} $ of ${\bfm g}$. If ${\bfm g}\in\mathds{C}(\Omega^+)$  takes on the form $\rho(x_1){\bfm \varphi}(\omega)$ for some $\rho\in
C^\infty_c([0;+\infty)) $ and ${\bfm \varphi}\in \mathcal{C}$, the entries of  $\partial {\bfm g} $ are given by
$$\partial_1 {\bfm g}(x_1,\omega)= \partial_{x_1}{\bfm g}(x_1){\bfm \varphi}(\omega),\quad  \text{ and, for }i=2,\dots,d,\,\,\,\,\,\partial_i {\bfm g}(x_1,\omega)=\rho(x_1)D_i{\bfm \varphi}(\omega).$$

%- $\partial_1 {\bfm g}(x_1,\omega)$ is given by $\partial_{x_1}{\bfm g}(x_1){\bfm \varphi}(\omega)$,

%-  $\partial_i {\bfm g}(x_1,\omega)$, for $i=2,\dots,d$, coincides with the random partial derivative $D_i$ of the random function $\omega\in \Omega\mapsto 
%{\bfm g}(x_1,\omega)$, where $x_1$ is fixed. 

%\begin{equation}\label{differ}
%\begin{split}\partial_1 {\bfm g}(x_1,\omega)&=\lim_{h>0,h\to
%0}\frac{{\bfm g}(x_1+h,\omega)-{\bfm g}(x_1,\omega)}{h},\\
%\partial_i {\bfm g}(x_1,\omega)&=\lim_{h\to 0}\frac{{\bfm
%g}(x_1,\tau_{h e_i}\omega)-{\bfm g}(x_1,\omega)}{h},\quad i=2,\dots,d.
%\end{split}
%\end{equation}
%In other words, if we fix $\omega\in \Omega$, $\partial_1 {\bfm g}$ is defined as the classical derivative of the function $x_1\in \R_+\mapsto {\bfm g}(x_1,\omega)$ whereas $\partial_i {\bfm g}$ (for $i=2,\dots,d$) coincides with the random partial derivative $D_i$ of the random function $\omega\in \Omega\mapsto 
%{\bfm g}(x_1,\omega)$, where $x_1$ is fixed. 

We define on $\mathds{C}(\Omega^+)$ the norm
\begin{equation}\label{norm}
  N({\bfm g})^2=|{\bfm g}(0,\cdot)|_2^2+\int_{\Omega^+}|\partial
{\bfm g}|_2^2\,d\mu^+,
\end{equation}
which is a sort of Sobolev norm on $\Omega^+$, and $\mathds{W}^1$ as the closure of $\mathds{C}(\Omega^+)$ with respect
to the norm $N$ ($\mathds{W}^1$ is thus an analog of Sobolev spaces on $ \Omega^+$). Obviously, the mapping $$P:\mathds{W}^1\ni{\bfm
g}\mapsto {\bfm g}(0,\cdot)\in L^2(\Omega) $$ is continuous (with
norm equal to 1) and stands, in a way, for the trace operator on $\Omega^+ $. Equip the topological dual space
$(\mathds{W}^1)'$ of $\mathds{W}^1$ with the dual norm $N'$. The adjoint $P^*$ of $P$ is
given by $P^*:{\bfm \varphi}\in L^2(\Omega)\mapsto P^*({\bfm
\varphi})\in (\mathds{W}^1)'$ where the mapping $P^*({\bfm \varphi})$ exactly matches $$P^*({\bfm \varphi}):{\bfm g}\in
\mathds{W}^1\mapsto ({\bfm g},P^*{\bfm \varphi})=({\bfm \varphi},{\bfm g}(0,\cdot) )_2 .$$
To sum up, we have constructed a space of test functions $\mathds{C}(\Omega^+)$, which is dense in $\mathds{W}^1$ for the norm $N$, and a trace operator on $\mathds{W}^1$. 

We further stress that a function ${\bfm g}\in \mathds{W}^1 $ satisfies  $\partial {\bfm g}=0 $ if and only if   we have ${\bfm g}(x_1,\omega)={\bfm f}(\omega)$ on $\Omega^+$ for some function ${\bfm f}\in L^2 (\Omega)$ invariant under the translations $\{\tau_x; x\in\{0\}\times\R^{d-1}\}$. For that reason, we introduce the $\sigma$-field ${\cal G}^*\subset {\cal G}$ generated by the subsets of $\Omega$ that are invariant under the translations $\{\tau_x; x\in\{0\}\times\R^{d-1}\}$, and the
conditional expectation $\M_1$ with respect to ${\cal G}^*$.

We now focus on the construction of the symmetric bilinear form and the resolvent family associated to $H_\gamma$. For each random function ${\bfm \varphi}$ defined on $\Omega$, we associate a function ${\bfm \varphi}^+$ defined on $\Omega^+$ by
$$\forall(x_1,\omega)\in \Omega^+,\quad  {\bfm \varphi}^+(x_1,\omega)={\bfm \varphi}(\tau_{x_1}\omega).$$ Hence, we can associate to the random matrix ${\bfm a}$ %and $\bar{{\bfm  H}}$ 
(defined in Section \ref{sec:intro}% and \ref{notations}
) the corresponding matrix-valued function ${\bfm a}^+$ %and $\bar{{\bfm H}}^+$ 
defined on $\Omega^+$. 
Then, for any $\lambda>0$, we define on $\mathds{W}^1\times
\mathds{W}^1$ the following symmetric bilinear form
\begin{align}\label{defbl1}
B_\lambda({\bfm g},{\bfm h}) & =\lambda(P{\bfm g},P{\bfm
h})_2+\frac{1}{2}\int_{\Omega^+}{\bfm a}^+%+\bar{{\bfm H}}^+)
_{ij}\,\partial_i{\bfm g}\,\partial_j
{\bfm
h}\,d\mu^+.%\\\label{defbl2}
%B_\lambda^*({\bfm g},{\bfm h}) & =\lambda(P{\bfm g},P{\bfm
%h})_2+\frac{1}{2}\int_{\Omega^+}({\bfm a}^+-\bar{{\bfm
%H}}^+)_{ij}\,\partial_i {\bfm g}\,\partial_j
%{\bfm h}\,d\mu^+
\end{align}
From \eqref{unif_ellip}, % the boundedness and the antisymmetry of $\bar{{\bfm H}}$,
 it is readily seen that it is continuous and coercive on
$\mathds{W}^1\times \mathds{W}^1$.  From the Lax-Milgram theorem, it thus defines a continuous resolvent family
$G_\lambda:(\mathds{W}^1)'\to \mathds{W}^1$ %and $G^*_\lambda:(\mathds{W}^1)'\to \mathds{W}^1$ 
such that:
\begin{equation}\label{resolG}
 \forall F\in (\mathds{W}^1)', \,\forall {\bfm g}\in \mathds{W}^1,
\quad B_\lambda(G_\lambda F,{\bfm g})=({\bfm g},F).%\quad \text{ and }\quad  B^*_\lambda(G^*_\lambda F,{\bfm g})=F({\bfm g}).
\end{equation} %Note that the restriction $G_\lambda:\mathds{W}^1\to \mathds{W}^1$ is self-adjoint.
For each
$\lambda>0$, we then define
the operator 
\begin{equation}\label{defrl}
R_\lambda:\begin{array}[t]{ccc}
L^2(\Omega)  &  \to &  L^2(\Omega) \\
{\bfm \varphi} & \mapsto & P G_\lambda  P^*({\bfm \varphi})
\end{array}.%\quad \text{ and }\quad R_\lambda^*:\begin{array}[t]{ccc}
%L^2(\Omega)  &  \to &  L^2(\Omega) \\
%{\bfm \varphi} & \mapsto & P G^*_\lambda  P^*({\bfm \varphi})
%\end{array}
\end{equation}
Given ${\bfm \varphi}\in L^2(\Omega)$, we can plug ${\bfm F} =P^*{\bfm \varphi}$ into \eqref{resolG}  and we get
\begin{equation}\label{resolR}
\forall {\bfm g}\in \mathds{W}^1,
\quad B_\lambda(G_\lambda P^*{\bfm \varphi},{\bfm g})=({\bfm g},P^*{\bfm \varphi}),%\quad \text{ and }\quad  B^*_\lambda(G^*_\lambda P^*{\bfm \varphi},{\bfm g})=({\bfm g},P^*{\bfm \varphi}),
\end{equation}
that is, by using \eqref{defbl1}:% and \eqref{defbl2}: 
\begin{align}
\label{resolR1}\forall {\bfm g}\in \mathds{W}^1,\quad \lambda(R_\lambda{\bfm \varphi},P{\bfm
g})_2+\frac{1}{2}\int_{\Omega^+}{\bfm a}^+%+\bar{{\bfm H}}^+)
_{ij}\,\partial_i(G_\lambda P^*{\bfm \varphi})\,\partial_j
{\bfm
g}\,d\mu^+=({\bfm g}(0,\cdot),{\bfm \varphi})_2,
%\\\label{resolR2}\forall {\bfm g}\in \mathds{W}^1,\quad \lambda(R_\lambda^*{\bfm \varphi},P{\bfm
%g})_2+\frac{1}{2}\int_{\Omega^+}({\bfm a}^+-\bar{{\bfm
%H}}^+)_{ij}\,\partial_i (G_\lambda^* P^*{\bfm \varphi})\,\partial_j
%{\bfm g}\,d\mu^+=({\bfm g},P^*{\bfm \varphi}).
\end{align}

The following proposition summarizes the main properties of the operators $(R_\lambda)_{\lambda>0}$, and in particular their ergodic properties:
\begin{proposition}\label{erglocal}
The family $(R_\lambda)_\lambda$ %and $(R_\lambda^*)_\lambda$ are both 
is a strongly continuous resolvent family, and:\\
1) the operator $R_\lambda$ is self-adjoint.\\
2) given ${\bfm
\varphi}\in L^2(\Omega)$ and $\lambda>0$, we have: $${\bfm
\varphi}\in {\rm Ker}(\lambda R_\lambda -{\rm
I})\,\,\,\Leftrightarrow\,\,\,{\bfm \varphi}=\M_1[{\bfm \varphi}]%\,\,\,\Leftrightarrow\,\,\,{\bfm
%\varphi}\in {\rm Ker}(\lambda R_\lambda -{\rm
%I})
.$$
3) for each function ${\bfm \varphi}\in L^2(\Omega)$, $|\lambda
R_\lambda{\bfm \varphi}-\M_1[{\bfm \varphi}] |_2\to 0$ as $\lambda\to 0$.
\end{proposition}

\vspace{2mm}
The remaining part of this section is concerned with the regularity properties of $G_\lambda P^*{\bfm
\varphi}$. 

\begin{proposition}\label{propregu}
Given ${\bfm \varphi}\in {\cal C}$, the trajectories of $G_\lambda P^*{\bfm \varphi}$ are smooth. More precisely, we can
find ${\bfm N}\subset \Omega $ satisfying $\mu({\bfm N})=0 $ and such
that $\forall \omega\in \Omega\setminus {\bfm N} $, the function
$$\tilde{u}_\omega:x=(x_1,y)\in \bar{D}\mapsto G_\lambda P^*{\bfm \varphi}(x_1,\tau_{(0,y)}\omega)$$ belongs to
$C^\infty(\bar{D}) $. Furthermore, it is a classical solution to the
problem:
\begin{equation}\label{prob:halfphi}
 \left\{\begin{array}{l}
L^\omega \tilde{u}_\omega(x)=0, \, x\in D,
  \\ \lambda \tilde{u}_\omega(x)-{\bfm \gamma}_i(\tau_x\omega)\partial_{x_i}\tilde{u}_\omega(x)={\bfm \varphi}(\tau_x\omega)
  ,\,x\in \partial D.
    \end{array}\right.
\end{equation}
\end{proposition}

In particular, the above proposition proves that $(R_\lambda)_\lambda$ is the resolvent family associated to the operator $H_\gamma$. This family also satisfies the maximum principle:
\begin{proposition}{\rm (Maximum principle).}\label{max}
Given ${\bfm \varphi}\in {\cal C}$ and $\lambda>0$, we have: $$ |G_\lambda P^*{\bfm \varphi}|_{L^\infty(\Omega^+)}\leq \lambda^{-1}|{\bfm \varphi}|_\infty.$$
\end{proposition}
%Our last result is concerned with the maximum principle for the
%function ${\bfm u}=G_\lambda P^*{\bfm \varphi} $.
%
%\begin{proposition}
%For any ${\bfm \varphi}\in {\cal C}$, we have $|\lambda R_\lambda
%{\bfm \varphi}|_\infty\leq |{\bfm \varphi}|_\infty $.
%\end{proposition}
%
%\noindent \textit{Proof.} We follow the proof of the classical maximum
%principle and choose a function $G\in C^\infty(\R)$ such that
%$|G'(s)|\leq M $ $\forall s\in\R$, $G$ is strictly increasing on
%$(0,+\infty)$ and $G(s)=0$ $\forall s \leq 0 $. We define $K>|{\bfm
%\varphi}|_\infty/\lambda$ and consider ${\bfm u}=G_\lambda P^*{\bfm
%\varphi}$. It is readily seen that ${\bfm v}=G({\bfm u}-K)\in
%\mathds{W}^1$ and that $\partial {\bfm v}=\partial{\bfm u} G'({\bfm
%u})$. Plugging ${\bfm g}={\bfm v} $ into the resolvent equation
%$\forall {\bfm g}\in \mathds{W}^1$, $B_\lambda({\bfm u},{\bfm
%g})=({\bfm \varphi},{\bfm g})_2 $ leads to
%$$\lambda ({\bfm u}(0,\cdot),G({\bfm u}(0,\cdot)-K))_2+\frac{1}{2}\int_{\R_+}\M\big[a_{ij}(x_1,\omega)
%\partial_i{\bfm u}\partial_j{\bfm u}G'({\bfm u}-K)\big]\,dx_1=({\bfm \varphi},G({\bfm u}(0,\cdot)-K))_2 ,$$
%in such a way that $ (\lambda {\bfm u}(0,\cdot)-{\bfm
%\varphi},G({\bfm u}(0,\cdot)-K))_2\leq 0$. Hence $ {\bfm
%u}(0,\cdot)\leq K$. Since $K>|{\bfm \varphi}|_\infty/\lambda$ is
%chosen arbitrarily, we deduce $ \mu$ a.s. $\lambda R_\lambda {\bfm
%\varphi}\leq |{\bfm \varphi}|_\infty$. We complete the proof by
%applying the same argument to the function $-{\bfm \varphi}$. \qed

%%%%%%%%%%%%%%%%%%%%%%%%%%%%%%%%%%%%%%%%%%%%%%%%%%%%%%%%%%%%%%%%%%%%%%%%%%%%%%%%%%%%%%%%%%%%%%
\subsection{Ergodic theorems}\label{ergth} 

As already explained, the ergodic problems that we have solved in the previous section lead to establishing  ergodic theorems for the process $X^\epsilon$. The strategy of the proof is the following. First we work under $\bar{\P}^{\varepsilon*}$ to use the existence of the IPM (see Section \ref{ipm}). By adapting a classical scheme, we derive from Propositions \ref{ergprob} and \ref{erglocal}   ergodic theorems under $\bar{\P}^{\varepsilon*}$ both for  the process $X^\epsilon$ and for the local time $K^\epsilon$:

\begin{theorem}\label{theorem_ergodic}
For each function ${\bfm f}\in L^1(\Omega)$ and $T>0$, we have
\begin{equation}\label{eq_ergodictheorem}
\lim_{\varepsilon\to 0}\bar{\E}^{\varepsilon *}\Big[\sup_{0 \leq t \leq
T}\big|\int_0^t{\bfm f}(\tau_{X^{\varepsilon}_r/\varepsilon}\omega)\,dr-t\M[{\bfm f}]\big|\Big]=0.
\end{equation}
\end{theorem}

\begin{theorem}\label{theorem_erglocal}
If ${\bfm f}\in L^2(\Omega) $, the following convergence holds
\begin{equation}\label{cvlocal}
  \lim_{\varepsilon\to 0}\bar{\E}^{\varepsilon *}\Big[\sup_{0\leq t \leq T}\big|\int_0^t{\bfm f}(\tau_{X^\varepsilon_r/\varepsilon}\omega)\,dK^\varepsilon_r
  -\M_1[{\bfm f}](\omega)K^\varepsilon_t\big|\Big]=0.
\end{equation}
\end{theorem}

%\begin{remark*} 
%In the proof of Theorem \ref{theorem_erglocal} and unlike   classical situations to establish ergodic theorems, we point out that we are faced with the lack of stationarity of the correctors (here the function $G_\lambda P^*{\bfm f} $) along the $x_1$-direction. The lack of stationarity of the correctors make them grow faster than the usual sublinear growth. This is not something that can be changed: whatever the notion of solution you give to problem \eqref{uphi}, you will face such a situation. To overcome that difficulty, we make sure that the decay of the IPM at infinity can compensate for that unusual growth. This is the reason why we impose the function $V$ to be of the type \eqref{defV}: it ensures a sufficient decay of the IPM at infinity (this is the meaning of  \eqref{V3} in the proof of Theorem \ref{theorem_erglocal}).
%\end{remark*}

Finally we deduce that the above theorems remain valid under  $\bar{\P}^{\varepsilon}$ thanks to \eqref{girsanov}.

\begin{theorem}\label{main_erg}
1) Let $({\bfm f}_\varepsilon)_\varepsilon$ be a family converging towards ${\bfm f}$ in $L^1(\Omega)$. For each fixed $ \delta>0$ and $T>0$, the following convergence holds
\begin{equation}\label{main_erg1}
\lim_{\varepsilon\to 0}\bar{\P}^{\varepsilon }\Big[\sup_{0 \leq t \leq
T}|\int_0^t{\bfm f}_\varepsilon(\tau_{X^{\varepsilon}_r/\varepsilon}\omega)\,dr-t\M[{\bfm f}]|\geq \delta\Big]=0.
\end{equation}

2)  Let $({\bfm f}_\varepsilon)_\varepsilon$ be a family converging towards ${\bfm f}$ in $L^2(\Omega)$. For each fixed $ \delta>0$ and $T>0$, the following convergence holds
\begin{equation}\label{main_erg2}
  \lim_{\varepsilon\to 0}\bar{\P}^{\varepsilon }\Big[\sup_{0\leq t \leq T}\big|\int_0^t{\bfm f}_\varepsilon(\tau_{X^\varepsilon_r/\varepsilon}\omega)\,dK^\varepsilon_r
  -\M_1[{\bfm f}]K^\varepsilon_t\big|\geq \delta\Big]=0.
\end{equation}
\end{theorem}

%%%%%%%%%%%%%%%%%%%%%%%%%%%%%%%%%%%%%%%%%%%%%%%%%%%
\subsection{Construction of the correctors}\label{corr}

Even though we have established ergodic theorems, this is not enough to find the limit of equation \eqref{SDEintro} because of the highly oscillating term $\varepsilon^{-1}{\bfm b}(\tau_{X^\varepsilon_t/\varepsilon}\omega)\,dt $. To get rid of this term, the ideal situation is to find a stationary solution ${\bfm u}^i:\Omega\to \R$ to the equation
 \begin{equation}\label{eq:resolven}
 -{\bfm L}{\bfm u}^i={\bfm b}_i.
\end{equation}
Then, by applying the It\^o formula to the function ${\bfm u}^i $, it is readily seen that the contribution of the term $\varepsilon^{-1}{\bfm b}_i(\tau_{X^\varepsilon_t/\varepsilon}\omega)\,dt $ formally reduces to a stochastic integral and a functional of the local time, the limits of which can be handled with the ergodic theorems \ref{main_erg}.

The problem is that the lack of compactness of a random medium makes you cannot find a stationary solution to \eqref{eq:resolven}. As already suggested in the literature, a good approach is to add some coercivity to the problem \eqref{eq:resolven} and define, for $i=1,\dots,d$ and $\lambda>0$, the solution ${\bfm u}_\lambda^i$ of the resolvent equation
\begin{equation}\label{eq:resolvent}
  \lambda {\bfm u}^i_\lambda-{\bfm L}{\bfm u}^i_\lambda={\bfm b}_i.
\end{equation}
If we let $\lambda$ go to $0$ in of \eqref{eq:resolvent}, the solution ${\bfm u}^i_\lambda$ should provide a good approximation of the solution of \eqref{eq:resolven}. Actually, it is hopeless to prove the convergence of the family $({\bfm u}^i_\lambda)_\lambda$ in some $L^p(\Omega)$-space because, in great generality, there is no stationary $L^p(\Omega)$-solution to \eqref{eq:resolven}. However we can prove the convergence towards $0$ of the term $ \lambda {\bfm u}^i_\lambda$ and the convergence of the gradients $D  {\bfm u}^i_\lambda $:
\begin{proposition}\label{prop:correctors}
There exists ${\bfm \zeta}^i\in (L^2(\Omega))^d$ such that
\begin{equation}\label{conv:xhi}
  \lambda|{\bfm u}^i_\lambda|_2^2+|D{\bfm
  u}^i_\lambda-{\bfm \zeta}^i|_2\rightarrow 0,\quad \text{ as }\lambda \to 0.
\end{equation}
\end{proposition}
As we shall see in Section \ref{homo}, the above convergence is enough to carry out the homogenization procedure. The functions ${\bfm \zeta}^i$ ($i\leq d$) are involved in the expression of the coefficients of the homogenized equation \eqref{SDElimit}. For that reason, we give some further qualitative description of these coefficients:

\begin{proposition}\label{prop:correctors2}
Define the random matrix-valued function ${\bfm \zeta}\in
L^2(\Omega;\R^{d\times d})$ by its entries ${\bfm \zeta}_{ij}={\bfm \zeta}_i^j =\lim_{\lambda\to 0}D_i{\bfm
  u}^j_\lambda$. Define the matrix  $\bar{A}$ and the d-dimensional vector $\bar{\Gamma} $ by 
 \begin{align}\label{derbarA} 
 \bar{A} & =\M[(\mathrm{I}+{\bfm \zeta}^*){\bfm a}(\mathrm{I}+{\bfm \zeta})],\,\,\text{which also matches }\M[(\mathrm{I}+{\bfm \zeta}^*){\bfm a}],
\\
\label{derbarG}
   \bar{\Gamma} & =\M[({\rm I}+{\bfm
\zeta}^*){\bfm \gamma}]\in\R^d ,
\end{align}
where ${\rm I}$ denotes the $d$-dimensional identity matrix. Then $\bar{A}$ obeys the variational formula:
  \begin{equation}\label{varform}
\forall X\in\R^d,\quad   X^* \bar{A}X=\inf_{{\bfm \varphi}\in \mathcal{C}}\M[(X+D{\bfm
\varphi})^*{\bfm a}(X+D{\bfm\varphi})].
\end{equation}
Moreover, we have $ \bar{A}\geq \Lambda\mathrm{I} $ (in the sense of symmetric nonnegative matrices) and the first component $\bar{\Gamma}_1$ of $\bar{\Gamma}$ satisfies $\bar{\Gamma}_1\geq \Lambda$. Finally, $\bar{\Gamma}$ coincides with the orthogonal projection
$\M_1[({\rm I}+{\bfm \zeta}^*){\bfm \gamma}] $.
\end{proposition}

In particular, we have established that the limiting equation \eqref{SDElimit} is not degenerate, namely that the diffusion coefficient $\bar{A}$ is invertible and that the pushing of the reflection term $\bar{\Gamma}$ along the normal to $\partial D$ does not vanish.

%%%%%%%%%%%%%%%%%%%%%%%%%%%%%%%%%%%%%%%%%%%%%%%%%%%
\subsection{Homogenization}\label{homo}

Homogenizing \eqref{SDEintro} consists in proving that the couple of processes $(X^\epsilon,K^\epsilon)_\epsilon$ converges as $\epsilon\to 0$ (in the sense of Theorem \ref{mainth}) towards the couple of processes $(\bar{X},\bar{K})$ solution of the RSDE \eqref{SDElimit}.  We also remind the reader that, for the time being, we work with the function $\chi(x)=e^{-2V(x)}$. We shall see thereafter how the general case follows. 

First we show that the family $(X^\epsilon,K^\epsilon)_\epsilon$ is compact in some appropriate topological space. Let us introduce the space $D([0,T];\R_+)$ of nonnegative
right-continuous functions with left limits on $[0,T]$ equipped with
the S-topology of Jakubowski (see Appendix \ref{sec:jaku}). The space $C([0,T];\bar{D})$ is equipped with the sup-norm topology. We have:
\begin{proposition}\label{prop:tightness}
Under the law $\bar{\P}^{\varepsilon }$, the family of processes
$(X^\varepsilon)_\varepsilon$ is tight in $C([0,T];\bar{D})$, and the
family of processes $(K^\varepsilon)_\varepsilon$ is tight in
$D([0,T];\R_+)$.
\end{proposition}

The main idea of the above result is originally due to Varadhan and is exposed in \cite[Chap. 3]{olla} for stationary diffusions in random media. Roughly speaking, it combines exponential estimates for processes symmetric with respect to their IPM and the Garsia-Rodemich-Rumsey inequality. In our context, the pushing of the local time rises some further technical difficulties when the process $X^\epsilon$ evolves near the boundary. Briefly, our strategy to prove Proposition \ref{prop:tightness} consists in  applying the method \cite[Chap. 3]{olla} when the process $X^\epsilon$ evolves far from the boundary, say not closer to $\partial D$ than a fixed distance $ \theta$, to obtain a first class of tightness estimates. Obviously, these estimates depend on $\theta$. That dependence takes place in a penalty  term related to the constraint of evolving far from the boundary. Then we let $\theta$ go to $0$. The limit of the penalty term can be expressed in terms of the local time $K^\epsilon$ in such a way that we get tightness estimates for the whole process $X^\epsilon$ (wherever it evolves). Details are set out in the appendix \ref{app:tight}.

It then remains to identify each possible weak limit of the family $(X^\epsilon,K^\epsilon)_\epsilon$. To that purpose, we  introduce the corrector ${\bfm
u}_\lambda\in L^2(\Omega;\R^d)$, the entries of which are given, for $j=1,\dots,d$, by the solution ${\bfm
u}^{j}_\lambda$ to the resolvent equation $$ \lambda {\bfm
u}^{j}_\lambda-{\bfm L}{\bfm
u}^{j}_\lambda={\bfm b}_j.$$
 Let ${\bfm \zeta}\in L^2(\Omega;\R^{d\times d})$ be defined by ${\bfm \zeta}_{ij}=\lim_{\lambda\to 0}D_i{\bfm
u}^{j}_\lambda$ (see  Proposition \ref{prop:correctors}). As explained in  Section \ref{corr}, the function ${\bfm
u}_\lambda$ is used to get rid of the highly oscillating term $\varepsilon^{-1}{\bfm b}(\tau_{X^\varepsilon_t/\varepsilon}\omega)\,dt $ in \eqref{SDEintro} by appling the It\^o formula.
Indeed, since $\mu $-almost surely the function $\phi:x\mapsto {\bfm u}_\lambda(\tau_x\omega)$ satisfies
$\lambda \phi-L^\omega
\phi={\bfm b}(\tau_\cdot \omega )$ on $\R^d$,  the function
$x\mapsto u_\lambda(\tau_x\omega) $ is smooth (see \cite[Th. 6.17]{gilbarg}) and we can apply the It\^o formula to the function $x\mapsto \epsilon {\bfm
u}_\lambda(\tau_{x/\epsilon}\omega)$. We obtain
\begin{align}
\varepsilon d {\bfm u}_{\lambda}(\tau_{X^\varepsilon_t/\varepsilon}\omega)=&\frac{1}{\epsilon}{\bfm L}{\bfm u}_{\lambda}(\tau_{X^\varepsilon_t/\varepsilon}\omega)\,dt+D
{\bfm u}_{\lambda}^*{\bfm \gamma}(\tau_{X^\varepsilon_t/\varepsilon}\omega)\,dK^\varepsilon_t+D
{\bfm u}_{\lambda}^*{\bfm \sigma}(\tau_{X^\varepsilon_t/\varepsilon}\omega)\,dB_t\nonumber\\
\label{temp}=& \frac{1}{\epsilon}(\lambda{\bfm u}_{\lambda}-{\bfm b})(\tau_{X^\varepsilon_t/\varepsilon}\omega)\,dt+D
{\bfm u}_{\lambda}^*{\bfm \gamma}(\tau_{X^\varepsilon_t/\varepsilon}\omega)\,dK^\varepsilon_t+D
{\bfm u}_{\lambda}^*{\bfm \sigma}(\tau_{X^\varepsilon_t/\varepsilon}\omega)\,dB_t.
\end{align}
By summing the relations \eqref{temp} and \eqref{SDEintro} and by setting $\lambda=\epsilon^2$, we deduce:
\begin{align}\label{ito}
X^\varepsilon_t = &x-\varepsilon\big({\bfm u}_{\varepsilon^2}(\tau_{X^\varepsilon_t/\varepsilon}\omega)-{\bfm u}_{\varepsilon^2}(\tau_{X^\varepsilon_0/\varepsilon}\omega)\big)+\epsilon\int_0^t{\bfm u}_{\varepsilon^2}(\tau_{X^\varepsilon_r/\varepsilon}\omega)\,dr
\\&+\int_0^t({\rm I}+D
{\bfm u}_{\varepsilon^2}^*){\bfm \gamma}(\tau_{X^\varepsilon_r/\varepsilon}\omega)\,dK^\varepsilon_r+\int_0^t({\rm
I}+D
{\bfm u}_{\varepsilon^2}^*){\bfm \sigma}(\tau_{X^\varepsilon_r/\varepsilon}\omega)\,dB_r.\nonumber\\
&\equiv x-
G^{1,\varepsilon}_t+G^{2,\varepsilon}_t+G^{3,\varepsilon}_t+M^\varepsilon_t.\nonumber
\end{align}
So we make the term $\varepsilon^{-1}{\bfm b}(\tau_{X^\varepsilon_t/\varepsilon}\omega)\,dt $ disappear 
at the price of modifying the stochastic integral and the integral with respect to the local time.  By using Theorem \ref{main_erg}, we should be able to identify their respective limits. The corrective terms $G^{1,\varepsilon} $ and $G^{2,\varepsilon}$ should reduce to $0$ as $\epsilon\to 0$. This is the purpose of the following proposition:

\begin{proposition}\label{m} For each subsequence of the family $(X^\varepsilon,K^\varepsilon)_\varepsilon$, we can extract a subsequence (still indexed with $\varepsilon>0$) such that:

1) under $\bar{\P}^\varepsilon$, the family of processes $(X^\varepsilon,M^\varepsilon,K^\varepsilon)_\varepsilon$ converges in law in $C([0,T]; \bar{D})\times C([0,T]; \R^d)\times D([0,T];\R_+)$ towards $(\bar{X},\bar{M},\bar{K})$, where $\bar{M}$ is a centered $d$-dimensional Brownian motion with covariance $$\bar{A}= \M[({\rm I}+{\bfm \zeta}^*){\bfm a}({\rm
I}+{\bfm \zeta})]$$ and $\bar{K}$ is a right-continuous increasing process.

2) the finite-dimensional distributions of the families $(G^{1,\varepsilon}_t)_\varepsilon$, $(G^{2,\varepsilon}_t)_\varepsilon$
and $(G^{3,\varepsilon}-\bar{\Gamma}K^\varepsilon)_\varepsilon$ converge towards $0$  in $\bar{\P}^\varepsilon$-probability, that is for each $t\in[0,T]$
$$\forall \delta>0,\quad \lim_{\varepsilon\to 0} \bar{\P}^\varepsilon\Big(|G^{i,\varepsilon}_t|>\delta\Big)=0 \,\,\,( i=1,2),\quad \lim_{\varepsilon\to 0}  \bar{\P}^\varepsilon\Big(|G^{3,\varepsilon}_t-\bar{\Gamma}K^\varepsilon_t|>\delta\Big)=0. $$
\end{proposition}

\noindent \textit{Proof.} 1) The tightness of $(X^\varepsilon,K^\varepsilon)$ results from Proposition \ref{prop:tightness}. To prove the tightness of the martingales $(M^\varepsilon)_\varepsilon$, it suffices to prove the tightness of the brackets $(<M^\varepsilon>)_\varepsilon$, which are given by  $$<M^\varepsilon>_t=\int_0^t
({\rm I}+D {\bfm u}_{\varepsilon^2}^*){\bfm a}({\rm I}+D
 {\bfm u}_{\varepsilon^2})(\tau_{X^\varepsilon_r/\varepsilon}\omega) \,dr.$$ Proposition \ref{prop:correctors} and Theorem
\ref{main_erg} lead to $<M^\varepsilon>_t\to \bar{A}t$ in probability in $C([0,T];\R^{d\times d}) $ where $\bar{A}=\M\big[({\rm I}+ {\bfm \zeta}^*){\bfm a}({\rm I}+
 {\bfm \zeta})\big]$. The martingales $(M^\varepsilon)_\varepsilon$ thus converge in law in $C([0,T];\R^{d}) $ towards a centered Brownian motion with covariance matrix $\bar{A}$ (see \cite{helland}). 

2) Let us investigate the convergence of $(G^{i,\varepsilon})_\varepsilon$ ($i=1,2$). From the Cauchy-Schwarz inequality, Lemma \ref{lem_invariant} and \eqref{conv:xhi}, we deduce:
$$\lim_{\varepsilon\to 0}\bar{\E}^{\varepsilon*}\Big[|\varepsilon{\bfm u}_{\varepsilon^2}(\tau_{X^\varepsilon_t/\varepsilon}\omega)|^2+|\int_0^t\varepsilon{\bfm u}_{\varepsilon^2}(\tau_{X^\varepsilon_r/\varepsilon}\omega)\,dr|^2\Big]\leq (1+t)\lim_{\varepsilon\to 0}( \varepsilon^2 |{\bfm
u}_{\varepsilon^2}|_2^2)=0. $$
We conclude with the help of \eqref{girsanov}.

Finally we prove the convergence of $(G^{3,\varepsilon})_\varepsilon$ with the help of Theorem \ref{main_erg}. Indeed, Proposition  \ref{prop:correctors} ensures the convergence of the family $(({\rm I}+D
{\bfm u}_{\varepsilon^2}^*){\bfm \gamma})_\epsilon$ towards $({\rm I}+{\bfm \zeta}^*){\bfm \gamma}$ in $L^2(\Omega)$ as $\epsilon\to 0 $. Furthermore we know from Proposition \ref{prop:correctors2} that $\bar{\Gamma}= \M_1[({\rm I}+{\bfm \zeta}^*){\bfm \gamma}] $. The convergence follows.\qed

\vspace{2mm}
Since the convergence of each term in \eqref{ito} is now established, it remains to identify the limiting equation. From Theorem \ref{thjaku}, we can find a countable subset $\mathcal{S}\subset[0,T[$ such that the finite-dimensional distributions of the process $(X^\varepsilon,M^\varepsilon,K^\varepsilon)_\varepsilon$ converge along $[0,T]\setminus  \mathcal{S}$. So we can pass to the limit in \eqref{ito} along $s,t\in [0,T]\setminus {\cal S}$ ($s<t$), and this  leads to
\begin{equation}\label{EDSlimit}
\bar{X}_t=\bar{X}_s+\bar{A}^{1/2}(\bar{B}_t-\bar{B}_s)+\bar{\Gamma}(\bar{K}_t-\bar{K}_s).
\end{equation}
Since \eqref{EDSlimit} is valid for $s,t\in[0,T]\setminus {\cal S}$ (note that this set is dense and contains $T$) and since the processes are at least right continuous,
\eqref{EDSlimit} remains valid on the whole interval $[0,T]$. As a
by-product, $\bar{K}$ is continuous and the convergence of
$(X^\varepsilon,M^\varepsilon,K^\varepsilon)_\varepsilon$  actually holds in the
space $C([0,T]; \bar{D})\times C([0,T]; \R^d)\times C([0,T];\R_+)$ (see Lemma \ref{vc}).

It remains to prove that $\bar{K}$ is associated to $\bar{X}$ in the
sense of the Skorokhod problem, that is to establish that  $\{\text{Points of increase of }\bar{K}\}\subset \{t;\bar{X}^1_t=0\}$ or $\int_0^T\bar{X}^1_r\,d\bar{K}_r=0$. This results from the fact that $\forall\varepsilon>0 $\, $\int_0^TX^{1,\varepsilon}_r\,dK^\varepsilon_r=0$ and Lemma \ref{contmap}.
Since uniqueness in law holds for the solution $(\bar{X},\bar{K})$ of
Equation \eqref{EDSlimit} (see \cite{watanabe}), we have proved that
each converging subsequence of the family $(X^\varepsilon,K^\varepsilon)_\varepsilon$
converges in law in $C([0,T];\bar{D}\times\R_+)$ as $\varepsilon\to
0 $ towards the same limit (the unique solution $(\bar{X},\bar{K}) $ of
\eqref{SDElimit}). As a consequence, under $\bar{\P}^\varepsilon$, the whole sequence $(X^\varepsilon,K^\varepsilon)_\varepsilon$ converges in law towards the couple  $(\bar{X},\bar{K})$ solution of \eqref{SDElimit}.

%The proof of Corollary \ref{coromain}  results from the fact that the unique viscosity solution of \eqref{EDP} has the probabilistic interpretation (see \cite{pardouxzhang})
%$$u_\varepsilon(t,x)=\E^\varepsilon_x\big[e^{-\lambda K^\varepsilon_t}f(X^\varepsilon_t)-\int_0^tg(X^\varepsilon_r)e^{-\lambda K^\varepsilon_r}\,dK^\varepsilon_r\big].$$
% Theorem \ref{mainth} and Proposition \ref{propg} (note that the parameter $\lambda$ makes the term $\int_0^tg(X^\varepsilon_r)e^{-\lambda K^\varepsilon_r}\,dK^\varepsilon_r$ bounded by $\lambda^{-1}\sup_{\partial D}|g|$ ).

\subsubsection*{Replication method}
%%%%%%%%%%%%%%%%%%%%%%%%%%%%%%%%%%%%%%
Let us use the shorthands $C_D$ and $C_+$ to denote the spaces $C([0,T],\bar{D})$ and $C([0,T],\R_+)$ respectively.  Let $\bar{\E}$ denote the expectation with respect to the law $\bar{\P}$ of the process $(\bar{X},\bar{K}) $ solving the RSDE \eqref{SDElimit} with initial distribution $\bar{\P}(\bar{X}_0\in dx)=e^{-2V(x)}dx$. From \cite{watanabe}, the law $\bar{\P}$ coincides with the averaged law $\int_{\bar{D}}\bar{\P}_x(\cdot)e^{-2V(x)}dx$ where $\bar{\P}_x$ denotes the law of $(\bar{X},\bar{K})$ solving \eqref{EDSlimit} and starting from $x\in\bar{D}$.

We sum up the results obtained previously. We have proved the convergence, as $\varepsilon\to 0$, of $\bar{\E}^\varepsilon[F(X^\varepsilon,K^\varepsilon)]$ towards $\bar{\E}[F(\bar{X},\bar{K})] $, for each continuous bounded function $F:C_D\times C_+\to \R$. This convergence result is often called annealed because $\bar{\E}^\varepsilon$ is the averaging of the law $\P^\epsilon_x $ with respect to the probability measure $\P^*_D$.

In the classical framework of Brownian motion driven SDE in random media (i.e. without reflection term in \eqref{SDEintro}), it is plain to see that the annealed convergence of $X^\varepsilon$  towards a Brownian motion implies that, in $\P^*_D$-probability, the law $ \P^\varepsilon_x$ of $X^\varepsilon$ converges towards that of a Brownian motion. To put it simply, we can drop the averaging with respect to $\P^*_D$ to obtain a convergence in probability, which is a stronger result. Indeed, the convergence in law towards $0$ of the correctors (by analogy, the terms $ G^{1,\varepsilon},G^{2,\varepsilon}$ in \eqref{ito}) implies their convergence in probability towards $0$. Moreover the convergence in $\P^*_D$-probability of the law of the martingale term $M^\varepsilon$ in \eqref{ito} is obvious since we can apply \cite{helland} for $\P^*_D$-almost every $(x,\omega)\in\bar{D}\times\Omega$. 

In our case, the additional term $G^{3,\varepsilon}$ puts an end to that simplicity: this term converges, under the annealed law $\bar{\P}^\varepsilon$, towards a random variable $\bar{\Gamma}\bar{K} $, but there is no obvious way to switch annealed convergence for convergence in probability. That is the purpose of the computations below. 
\begin{remarkopen}The above remark also raises the open problem of proving a so-called quenched homogenization result, that is to prove the convergence of $X^\epsilon$ towards a reflected Brownian motion for almost every realization $\omega$ of the environment and every starting point $x\in\bar{D}$. The same arguments as above show that a quenched result should be much more difficult than in the stationary case \cite{sznitman}.\end{remarkopen}

So we have to establish the convergence in $\P^*_D $-probability of $\E^\varepsilon_x[F(X^\varepsilon,K^\varepsilon)]$ towards $\bar{\E}_x [F(\bar{X},\bar{K})] $ for each continuous bounded function $F:C_D\times C_+\to \R$. Obviously, it is enough to prove the convergence of $\E^\varepsilon_x[F(X^\varepsilon,K^\varepsilon)]$ towards $\bar{\E}_x [F(\bar{X},\bar{K})] $ in $L^2(\bar{D}\times \Omega, \P^*_D)$. By using a specific feature of Hilbert spaces, the convergence is established if we can prove the convergence of the norms 
\begin{equation}\label{cvcarre}
\M^*_D\Big[\big(\E^\varepsilon_x[F(X^\varepsilon,K^\varepsilon)]\big )^2\Big]\to\M^*_D\Big[\big(\bar{\E}_x [F(\bar{X},\bar{K})]\big )^2\Big]\quad \text{as }\varepsilon\to 0,
\end{equation}
as well as the weak convergence. 
%\begin{equation}\label{cvcar}
%\lim_{\varepsilon\to 0}\M^*_D\Big[\big(\E^\varepsilon_x[F(X^\varepsilon,K^\varepsilon)]-\bar{\E}_x [F(\bar{X},\bar{K})]\big )^2\Big]=0
%\end{equation}
%The expression $\M^*_D\Big[\big(\E^\varepsilon_x[F(X^\varepsilon,K^\varepsilon)]-\bar{\E}_x [F(\bar{X},\bar{K})]\big )^2\Big]$ expands as $$\M^*_D\Big[\big(\E^\varepsilon_x[F(X^\varepsilon,K^\varepsilon)]\big )^2\Big]-2\M^*_D\Big[\E^\varepsilon_x[F(X^\varepsilon,K^\varepsilon)]\bar{\E}_x [F(\bar{X},\bar{K})]\Big]+\M^*_D\Big[\big(\bar{\E}_x [F(\bar{X},\bar{K})]\big )^2\Big].$$
%Clearly, the mapping $x\in \bar{D}\mapsto \bar{\E}_x \big[F(\bar{X},\bar{K})\big]$ is bounded and continuous so that
%$$G:(h,k)\in C_D\times C_+\mapsto F(h,k)\E_{h(0)} \big[F(\bar{X},\bar{K})\big]$$ is bounded and continuous. Hence the following convergence holds as $\varepsilon\to 0$
%$$\M^*_D\Big[\E^\varepsilon_x[F(X^\varepsilon,K^\varepsilon)]\bar{\E}_x [F(\bar{X},\bar{K})]\Big]= \bar{\E}^\varepsilon[G(X^\varepsilon,K^\varepsilon)]\to   \bar{\E} [G(\bar{X},\bar{K})]\Big]= \M^*_D\Big[(\bar{\E}_x \big[F(\bar{X},\bar{K})]\big)^2\Big]$$
%So, we just have to prove that 
%\begin{equation}\label{cvcarre}
%\M^*_D\Big[\big(\E^\varepsilon_x[F(X^\varepsilon,K^\varepsilon)]\big )^2\Big]\to\M^*_D\Big[\big(\bar{\E}_x [F(\bar{X},\bar{K})]\big )^2\Big]\quad \text{as }\varepsilon\to 0.
%\end{equation}
Actually we only need to establish \eqref{cvcarre} because the weak convergence results from Section \ref{homo} as soon as \eqref{cvcarre} is established. 

The following method is called replication technique because the above quadratic mean can be thought as of the mean of two independent copies of the couple $(X^\varepsilon,K^\varepsilon)$. We consider 2 independent Brownian motions $(B^1,B^2)$ and solve \eqref{SDEintro} for each Brownian motion. This provides two independant (with respect to the randomness of the Brownian motion) couples of processes $(X^{\varepsilon,1},K^{\varepsilon,1}) $ and $(X^{\varepsilon,2},K^{\varepsilon,2}) $.
Furthermore, we have $$\M^*_D\Big[\big(\E^\varepsilon_x[F(X^\varepsilon,K^\varepsilon)]\big )^2\Big]=\M^*_D\big[\E^\varepsilon_{xx}\big[F(X^{\varepsilon,1},K^{\varepsilon,1})F(X^{\varepsilon,2},K^{\varepsilon,2})\big]\big]$$ where $\E^\varepsilon_{xx}$ denotes the expectation with respect to the law 
$\P^\varepsilon_{xx}$ of the process $(X^{\varepsilon,1},K^{\varepsilon,1},X^{\varepsilon,2},K^{\varepsilon,2})$ when both $X^{\varepsilon,1}$ and $X^{\varepsilon,2}$ start from $x\in\bar{D}$. Under $\M^*_D\P^\varepsilon_{xx}$, the results of subsections \ref{ergo}, \ref{corr} and Proposition \ref{prop:tightness} remain valid since the marginal laws of each couple of processes coincide with $\bar{\P}_x^\varepsilon$. So we can repeat the arguments of subsection \ref{homo} and prove that the processes $(X^{\varepsilon,1},K^{\varepsilon,1},X^{\varepsilon,2},K^{\varepsilon,2})_\varepsilon$ converge in law in $C_D\times C_+\times C_D\times D_+$, under $\M^*_D\E^\varepsilon_{xx} $, towards a process $(\bar{X}^1,\bar{K}^1,\bar{X}^2,\bar{K}^2)$ satisfying:
\begin{equation}\label{indep}
\forall t\in[0,T],\quad \bar{X}^1_t=\bar{X}^1_0+A^{1/2}\bar{B}^1_t +\bar{\Gamma}\bar{K}^1_t,\quad \bar{X}^2_t=\bar{X}^2_0+A^{1/2}\bar{B}^2_t +\bar{\Gamma}\bar{K}^2_t,
\end{equation}
 where $(\bar{B}^1,\bar{B}^2)$ is a standard $2d$-dimensional Brownian motion and $\bar{K}^1,\bar{K}^2$ are the local times respectively associated to $\bar{X}^1,\bar{X}^2 $. Let $\bar{\P}$ denote the law of $(\bar{X}^1,\bar{K}^1,\bar{X}^2,\bar{K}^2)$ with initial distribution given by $\bar{P}(\bar{X}^1_0\in dx,\bar{X}^2_0\in dy)=\delta_x(dy)e^{-2V(x)}dx$ and $\bar{\P}_{xx}$ the law of $(\bar{X}^1,\bar{K}^1,\bar{X}^2,\bar{K}^2)$ solution of \eqref{indep} where both $\bar{X}^1$ and $\bar{X}^2$ start from $x\in\bar{D}$.
 To obtain \eqref{cvcarre}, it just remains to remark that 
\begin{align*}
\bar{\E}\big[F(\bar{X}^1,\bar{K}^1)F(\bar{X}^2,\bar{K}^2)\big] =&\int_{\bar{D}}\bar{\E}_{xx}\big[F(\bar{X}^1,\bar{K}^1)F(\bar{X}^2,\bar{K}^2)\big]e^{-2V(x)}dx\\
=&\int_{\bar{D}}\bar{\E}_{x}\big[F(\bar{X}^1,\bar{K}^1)\big]\bar{\E}_{x}\big[F(\bar{X}^2,\bar{K}^2)\big]e^{-2V(x)}dx,
\end{align*}
since, under $\bar{\P}_{xx}$, the couples $(\bar{X}^1,\bar{K}^1) $ and $(\bar{X}^2,\bar{K}^2)$ are adapted to the filtrations generated respectively by $\bar{B}^1$ and $\bar{B}^2$ and are therefore independent.\qed

% and $ \alpha=\M[{\bfm \sigma}{\bfm \sigma}^*]$
%\begin{equation}\label{qv}\int_0^t\left(\begin{array}{cc}
%A & A\nabla V(\bar{X}_r) \\ 
% \nabla V^*(\bar{X}_r)A&
% \nabla V^*(\bar{X}_r)\alpha  \nabla V(\bar{X}_r)\end{array}\right)\,dr,
% \end{equation}
%
%$$ \int_0^t\left(\begin{array}{cc}
%({\rm I}+D {\bfm u}_{\varepsilon^2}^*){\bfm a}({\rm I}+D
% {\bfm u}_{\varepsilon^2})^*(\tau_{X^\varepsilon_r/\varepsilon}\omega) &({\rm I}+D {\bfm u}_{\varepsilon^2}^*){\bfm a} (\tau_{X^\varepsilon_r/\varepsilon}\omega)\nabla V(X^\varepsilon_r) \\ 
% \nabla V^*(X^\varepsilon_r) {\bfm a}({\rm I}+D
% {\bfm u}_{\varepsilon^2})^*(\tau_{X^\varepsilon_r/\varepsilon}\omega)&
% \nabla V^*(X^\varepsilon_r) {\bfm a}(\tau_{X^\varepsilon_r/\varepsilon}\omega)  \nabla V(X^\varepsilon_r)\end{array}\right)\,dr.$$

%%%%%%%%%%%%%%%%%%%%%%%%%%%%%%%%%%%%%%%%%%%%%%%%%%%
\subsection{Conclusion}

We have proved Theorem \ref{mainth} for any function $\chi $ that can be rewritten as $\chi(x)=e^{-2V(x)} $, where $V:\bar{D}\to \R$ is defined in \eqref{defV}.
It is then plain to see that Theorem \ref{mainth} holds for any nonnegative function $\chi$ not greater than  $Ce^{-2V(x)}$, for some positive constant $C$ and some function $V$ of the type \eqref{defV}. %Indeed, each measurable subset $A\subset \Omega\times\bar{D}$ satisfies $\M\int_{\bar{D}}\one_A\chi(x)\,dx\leq C \M\int_{\bar{D}}\one_Ae^{-2V(x)}\,dx$. 
Theorem \ref{mainth} thus holds for any continuous function $\chi$ with compact support over $\bar{D}$. 

Consider now a generic function $\chi:\bar{D}\to \R_+$ satisfying $\int_{\bar{D}}\chi(x)\,dx=1$ and $\chi':\bar{D}\to \R_+$ with compact support in $\bar{D}$. Let $A^\varepsilon\subset \Omega\times \bar{D}$ be defined as  
$$ A^\varepsilon=\left\{(\omega,x)\in \Omega\times \bar{D};\big|\E^\varepsilon_x(F(X^\varepsilon,K^\varepsilon))-\E_x(F(\bar{X},\bar{K}))\big|\geq \delta\right\}.$$ 
From the relation $\M\int_{\bar{D}}\one_{A^\varepsilon}\chi(x)dx\leq \M\int_{\bar{D}}|\chi(x)-\chi'(x)|dx+ \M\int_{\bar{D}}\one_{A^\varepsilon}\chi'(x)dx$, we deduce
$$\limsup_{\epsilon\to 0}\M\int_{\bar{D}}\one_{A^\varepsilon}\chi(x)dx\leq \M\int_{\bar{D}}|\chi(x)-\chi'(x)|dx,$$ in such a way that the Theorem \ref{mainth} holds for $\chi$ by density arguments. The proof is completed.

%%%%%%%%%%%%%%%%%%%%%%%%%%%%%%%%%%%%%%%%%%%%%%%%%%%%%%%%%%%%%%%%%%%%%%%%%%%%%%%%%%%%%%%%%%%%%%%%%%%%%%%%%%%%%%%%%%%%%%%%%%%%%%%%%%%%%%%%%%%%%%%%%%%%%%%
\appendix
\section*{Appendix}
%%%%%%%%%%%%%%%%%%%%%%%%%%%%%%%%%%%%%%%%%%%%%%%%%%%%%%%%%%%%%%%%%%%%%%%%%%%%%%%%%%%%%%%%%%%%%%%%%%%%%%

%%%%%%%%%%%%%%%%%%%%%%%%%%%%%%%%%%%%%%%%%%%%%%%%%
\section{Preliminary results}\label{app:aux}
%%%%%%%%%%%%%%%%%%%%%%%%%%%%%%%%%%%%%%%%%%%%%%%%%
\textit{Notations: Classical spaces.} Given an open domain ${\cal O}\subset \R^n$ and $k\in \nat\cup\{\infty\}$, $C^k({\cal O})$ (resp. $C^k(\bar{{\cal O}})$, resp. $C^k_b(\bar{{\cal O}})$) denotes the space of functions admitting continuous derivatives up to order $k$ over ${\cal
O}$ (resp. over $\bar{\mathcal{O}}$, resp. with continuous bounded derivatives over $\bar{D}$). The spaces $C^k_c({\cal O})$ 
and $ C^k_c(\bar{{\cal O}})$ denote the subspaces of
$C^k(\bar{{\cal O}})$ whose functions respectively  have a compact support in
${\cal O}$ or have a compact support in $\bar{{\cal O}}$. Let $C^{1,2}_b$ denote the space of bounded functions $f:[0,T]\times\bar{D}\to\R$ admitting bounded and continuous derivatives $\partial_t f$, $\partial_xf$, $ \partial^2_{tx}f$ and $\partial^2_{xx}f$ on $[0,T]\times\bar{D}$.

\subsubsection*{Green's formula:}
We remind the reader of the Green formula (see \cite[eq.
6.5]{miranda}). %Given a smooth matrix-valued function
%$c:\bar{D}\to \R^{d\times d} $, and $(\varphi,\psi)\in
%C^2(\bar{D})\times C^1_c(\bar{D})$:
%\begin{equation}\label{green}
%\int_D\partial_{x_i}\big(c_{ij}\partial_{x_j}\varphi\big)\psi(x)\,dx=-\int_{\partial
%D}c_{1j}\partial_{x_j}\varphi\psi(x)\,dx-\int_Dc_{ij}\partial_{x_j}\varphi\partial_{x_i}\psi(x)\,dx.
%\end{equation}
We consider the following operator acting on $ C^2(\bar{D})$
\begin{equation}\label{gen:introV2}
{\cal
L}^\varepsilon_V=\frac{e^{2V(x)}}{2}\sum_{i,j=1}^d\partial_{x_i}\big(e^{-2V(x)}{\bfm a}_{ij}(\tau_{x/\varepsilon}\omega)\partial_{x_j}\,\big),
\end{equation}
where $V:\bar{D}\to \R$ is smooth. For any couple  $(\varphi,\psi)\in
C^2(\bar{D})\times C^1_c(\bar{D})$, we have 
\begin{align}\label{green2}
\int_D{\cal L}^{\varepsilon}_V\varphi (x)\psi(x)e^{-2V(x)}\,dx&
+\frac{1}{2}\int_D{\bfm a}
%+\bar{{\bfm H}})
_{ij}(\tau_{x/\varepsilon}\omega)\partial_{x_i}\varphi(x)\partial_{x_j}\psi(x)e^{-2V(x)}\,dx\nonumber\\&=
-\frac{1}{2}\int_{\partial D}{\bfm \gamma}_i(\tau_{x/\varepsilon}\omega)\partial_{x_i}
\varphi(x)\psi(x)e^{-2V(x)}\,dx.
\end{align}
Note that the Lebesgue
measure on ${\bar{D}}$ or $\partial D$ is indistinctly denoted by $dx$ since
the domain of integration avoids confusion.

\subsubsection*{PDE results:}

We also state some preliminary PDE results that we shall need in the forthcoming proofs:
\begin{lemma}\label{ladyzh}
For any  functions $f\in C^\infty_c(D)$ and $g,h\in C^\infty_b(\bar{D})$, there exists a unique classical solution $w_\varepsilon\in C^\infty([0,T];\bar{D})\cap C^{1,2}_b$ to the problem 
\begin{equation}\label{eqlady}
\partial_t w_\varepsilon={\cal
L}^\varepsilon_V w_\varepsilon +gw_\varepsilon+h\text{ on } [0,T]\times D,\quad {\bfm \gamma}_i(\tau_{\cdot/\varepsilon}\omega)\partial_{x_i}w_{\varepsilon}=0 \text{ on }
[0,T]\times\partial D,\quad  \text{ and }w_\varepsilon(0,\cdot)=f.
\end{equation}
\end{lemma}

\vspace{2mm}
\noindent \textit{Proof.} First of all, we remind the reader that all the coefficients involved in the operator ${\cal L}^\varepsilon_V$  belong to $C^\infty_b(\bar{D})$. From \cite[Th V.7.4]{lady}, we can find a unique generalized solution
$w'_{\varepsilon}$ in $C^{1,2}_b$ to the equation $$\partial_t w'_\varepsilon={\cal L}^\varepsilon_V
w'_\varepsilon+gw'_\varepsilon+{\cal L}^\varepsilon_V f+gf+h,\quad w'_{\varepsilon}(0,\cdot)=0 \text{ on }D,\quad {\bfm \gamma}(\tau_{\cdot/\varepsilon}\omega)\partial_{x_i}w'_{\varepsilon}=0 \text{ on }
[0,T]\times\partial D  .$$ From \cite[IV.§10]{lady}, we can
prove that $w'_\varepsilon $ is smooth up to the boundary. Then the function $$w_\varepsilon(t,x)=w'_\varepsilon(t,x)+f(x)\in C^\infty([0,T]\times
\bar{D})\cap C_b^{1,2}$$
% has $L^2(D)$-integrable derivatives up to order 2, and 
is a classical solution to the problem \eqref{eqlady}.
% Moreover $ w_\varepsilon$
%is bounded \cite[Th.V.7.3]{lady}.
% as well as its gradient \cite[Th.V.7.2]{lady}.
\qed

\begin{lemma}\label{FFK}
The solution $w_\varepsilon$ given by Lemma \ref{ladyzh} admits the following probabilistic representation: $\forall (t,x)\in [0,T]\times \bar{D},$
$$ w_\varepsilon(t,x)=\E^{\varepsilon*}_x\Big[f(X^\varepsilon_t)\exp\Big(\int_0^tg(X^\varepsilon_r)dr\Big)+\int_0^th(X^\varepsilon_r)\exp\Big(\int_0^rg(X^\varepsilon_u)du\Big)dr\Big]. $$
\end{lemma}

\vspace{2mm}
\noindent \textit{Proof.} The proof relies on the It\^o formula  (see for instance \cite[Ch. II, Th. 5.1]{ikeda} or \cite[Ch. 2, Th. 5.1]{freidlin}). It must be applied to the function $(r,x,y)\mapsto w_\varepsilon(t-r,x)\exp(y) $ and to the triple of processes $(r,X^\varepsilon_r,\int_0^rg(X^\varepsilon_u)du) $. Since it is a quite  classical exercise, we let the reader check the details.\qed

% Applying the It\^o formula to the function $(r,x,y)\mapsto w_\varepsilon(t-r,x)\exp(y) $ and to the triple of processes $(r,X^\varepsilon_r,\int_0^rg(X^\varepsilon_u)du) $ yields (see for instance \cite[Ch. II, Th. 5.1]{ikeda} or \cite[Ch. 2, Th. 5.1]{freidlin}):
%\begin{align*}
%d\Big(w_\varepsilon(t-r,X^{\varepsilon}_r)\exp\big(\int_0^rg(X^\varepsilon_u)du\big)\Big)=&\exp\big(\int_0^rg(X^\varepsilon_u)du\big)\big((-\partial_tw_\varepsilon+{\cal
%L}^\varepsilon_V w_\varepsilon)(t-r,X^\varepsilon_r)\big)dr\\
%&+\exp\big(\int_0^rg(X^\varepsilon_u)du\big)\partial_{x_i}w_\varepsilon(t-r,X^\varepsilon_r){\bfm \sigma}_{ij}(\tau_{X^\varepsilon_r/\varepsilon}\omega)dB^{*j}_r\\
%&+w_\varepsilon(t-r,X^\varepsilon)g(X^\varepsilon_r)\exp\big(\int_0^rg(X^\varepsilon_u)du\big)dr\\
%&+\exp\big(\int_0^rg(X^\varepsilon_u)du\big){\bfm \gamma}_i(\tau_{X^\varepsilon_r/\varepsilon}\omega)\partial_{x_i}w_\varepsilon(t-r,X^\varepsilon_r)dK^\varepsilon_r,
%\end{align*}
%that is, by using \eqref{eqlady},
%\begin{align*}
%f(X^{\varepsilon}_t)\exp\big(\int_0^tg(X^\varepsilon_u)du\big)=&w_\varepsilon(t,x)-\int_0^th(X^\varepsilon_r)\exp\big(\int_0^rg(X^\varepsilon_u)du\big)dr\\
%&+\int_0^t\exp\big(\int_0^rg(X^\varepsilon_u)du\big)\partial_{x_i}w_\varepsilon(t-r,X^\varepsilon_r){\bfm \sigma}_{ij}(\tau_{X^\varepsilon_r/\varepsilon}\omega)dB^{*j}_r.
%\end{align*}
%Since $w_\varepsilon\in C^{1,2}_b$, the stochastic integral is a martingale and its expectation reduces %to $0$. So it just remains to take the expectation in the above expression to prove the Lemma.\qed 

%%%%%%%%%%%%%%%%%%%%%%%%%%%%%%%%%%%%%%%%%%%%%%%%%
\section{Proofs of subsection \ref{ipm}}\label{app:ipm}
%%%%%%%%%%%%%%%%%%%%%%%%%%%%%%%%%%%%%%%%%%%%%%%%%

\noindent \textit{Proof of Lemma \ref{lem_invariant}.} 1) Fix $t>0$. First we suppose that we are given a deterministic function $f:\bar{D}\to \R$ belonging to $ C^\infty_c(D)$. From Lemma \ref{ladyzh}, there exists a classical bounded solution $w_\varepsilon\in C^\infty([0,t]\times
\bar{D})\cap C_b^{1,2}$  to the problem 
$$\partial_t w_\varepsilon={\cal
L}^\varepsilon_V w_\varepsilon \text{ on } [0,t]\times D,\quad {\bfm \gamma}_i(\tau_{\cdot/\varepsilon}\omega)\partial_{x_i}w_{\varepsilon}=0 \text{ on }
[0,t]\times\partial D,\quad  \text{ and }w_\varepsilon(0,\cdot)=f(\cdot),$$ where ${\cal
L}^\varepsilon_V $ is defined in \eqref{gen:introV2}. Moreover, Lemma \ref{FFK} provides the probabilistic representation: $$w_\varepsilon(t,x)=\E^{\varepsilon*}_x[f(X^{\varepsilon}_t)] .$$
%Obviously, the Feynman-Kac formula \cite[Ch.
%II, Th. 5.1]{freidlin} provides the probabilistic
%representation $w_\varepsilon(t,x)=\E^{\varepsilon*}_x[{\bfm
%f}(X^{\varepsilon}_t,\tau_{X^\varepsilon_t/\varepsilon}\omega )] $.
The Green formula \eqref{green2} then yields
\begin{align*}
\partial_t\int_Dw_\varepsilon(t,x) e^{-2V(x)}\,dx&= \int_D{\cal
L}^{\varepsilon}_Vw_\varepsilon(t,x) e^{-2V(x)}\,dx\\&=-\frac{1}{2} \int_{
\partial D}{\bfm \gamma}_i (\tau_{x/\varepsilon}\omega)\partial_{x_i}w_\varepsilon(t,x)
e^{-2V(x)}\,dx=0
\end{align*}
so that 
\begin{equation}\label{ipmtemp}
\int_{\bar{D}}\E^{\varepsilon*}_x[f(X^{\varepsilon}_t)]e^{-2V(x)}\,dx= \int_{\bar{D}}{\bfm f}(x)e^{-2V(x)}\,dx.
\end{equation}
It is readily seen that \eqref{ipmtemp} also holds if we only assume that $f$ is a bounded and continuous function over $\bar{D}$: it suffices to consider a sequence $(f_n)_n\subset C^\infty_c(D)$ converging pointwise towards $f$ over $D$. Since $f$ is bounded, we can assume that the sequence is uniformly bounded with respect to the sup-norm over $\bar{D}$. Since \eqref{ipmtemp} holds for $f_n$, it just remain to pass to the limit as $n\to \infty$ and apply the Lebesgue dominated convergence theorem. 

We have proved that the measure $e^{-2V(x)}\,dx$ is invariant for the Markov process $X^{\varepsilon}$ (under $\P^{\epsilon*}$). Its semi-group thus uniquely extends to a contraction semi-group on $L^1(\bar{D},e^{-2V(x)}\,dx)$.

Consider now ${\bfm f}\in L^1(\bar{D}\times \Omega;\P^*_D)$ and $\epsilon>0$. Then, $\mu$ almost surely, the mapping $x\mapsto {\bfm f}(x,\tau_{x/\epsilon}\omega)$ belongs to $L^1(\bar{D},e^{-2V(x)}\,dx)$. Applying  \eqref{ipmtemp} yields, $\mu$ almost surely, 
$$\int_{\bar{D}}\E^{\varepsilon*}_x[{\bfm
f}(X^{\varepsilon}_t,\tau_{X^\varepsilon_t/\varepsilon}\omega
)]e^{-2V(x)}\,dx= \int_{\bar{D}}{\bfm f}(x,\tau_{x/\varepsilon}\omega)e^{-2V(x)}\,dx. $$
It just remains to
integrate with respect to the measure $\mu$ and use the
invariance of $\mu$ under translations. 

Let us now focus on the second assertion. As previously, it suffices to establish
$$\int_{\bar{D}}\E^{\varepsilon *} _x\big[ \int_0^tf(X^\varepsilon_r)\,dK^\varepsilon_r\big]e^{-2V(x)}\,dx=t\int_{\partial D}f(x)e^{-2V(x)}\,dx$$
for some bounded continuous function $f:\partial D\to \R$. We can find a bounded continuous function $\tilde{f}:\bar{D}\to \R$ such that the restriction to $\partial D$ coincides with $f$ (choose for instance $\tilde{f}=f\circ p$ where $p:\bar{D}\to \partial D$ is the orthogonal projection along the first axis of coordinates).

Recall now that the local time
$K^\varepsilon_t$ is the density of occupation time at $\partial D$
(see \cite[Prop. 1.19]{cattiaux} with $\psi(x)=x_1$, $V_0={\bfm \gamma}$ and $a^2(x)=1$). Hence, by using \eqref{ipmtemp},
\begin{align*}
\int_{\bar{D}}\E^{\varepsilon *} _x\big[ \int_0^tf(X^\varepsilon_r)\,dK^\varepsilon_r\big]e^{-2V(x)}\,dx
= & \int_{\bar{D}}\E^{\varepsilon *} _x\big[\lim_{\delta\to 0}\delta^{-1}
\int_0^t\tilde{f}(X^\varepsilon_r)\mathds{1}_{[0,\delta]}(X^{1,\varepsilon}_r)\,dr\big]e^{-2V(x)}\,dx\\
= & \lim_{\delta\to 0} \int_{\bar{D}}\E^{\varepsilon *} _x\big[\delta^{-1}
\int_0^t \tilde{f}(X^\varepsilon_r)\mathds{1}_{[0,\delta]}(X^{1,\varepsilon}_r)\,dr\big]e^{-2V(x)}\,dx\\
= & t\lim_{\delta\to 0} \delta^{-1}\int_{\bar{D}}\tilde{f}(x)\mathds{1}_{[0,\delta]}(x_1)e^{-2V(x)}\,dx\\
=&t\int_{\partial D}f(x)e^{-2V(x)}\,dx.\qed
\end{align*}

\section{Proofs of subsection \ref{ergo}}\label{app:ergo}
%%%%%%%%%%%%%%%%%%%%%%%%%%%%%%%%%%%%%%%%%%%%%%%%%
\subsubsection*{Generator on the random medium associated to the
diffusion process inside $D$} 

\vspace{2mm}\noindent \textit{Proof of Proposition \ref{ergprob}.} The first statement is a particular case, for instance, of \cite[Lemma 6.2]{rhodes:06}. To follow the proof in \cite{rhodes:06}, omit the dependency on the parameter $y$, take ${\bfm H}=0$ and ${\bfm \Psi}={\bfm f}$. To prove the second statement, choose ${\bfm \varphi}={\bfm w}_\lambda$ in \eqref{wfe} and plug the relation
$$({\bfm f},{\bfm w}_\lambda)_2\leq |{\bfm f}|_2|{\bfm w}_\lambda|_2\leq 1/(2\lambda)|{\bfm f}|^2_2+(\lambda/2)|{\bfm w}_\lambda|_2^2 $$ into the right-hand side to obtain $\lambda |{\bfm w}_\lambda|_2^2 +\big({\bfm a}_{ij}D_i{\bfm w}_\lambda,D_j{\bfm w}_\lambda\big)_2\leq |{\bfm f}|^2_2/\lambda$. From \eqref{unif_ellip}, we deduce $\Lambda|D{\bfm w}_\lambda|_2^2\leq   |{\bfm f}|^2_2/\lambda$ and the result follows. \qed

\vspace{2mm}\noindent \textit{Proof of Lemma \ref{maxU}.} The proof is quite similar to that of Proposition \ref{max} below. So we let the reader check the details.\qed

\subsubsection*{Generator on the random medium associated to the
reflection term}  
%%%%%%%%%%%%%%%%%%%%%%%%%%%%%%
\noindent \textit{Proof of Proposition \eqref{erglocal}.} The resolvent properties of the family $(R_\lambda)_\lambda$ are readily derived from those of  the family
$(G_\lambda)_\lambda$. % and $G_\lambda^*$. 

So we first prove 1). Consider ${\bfm \varphi},{\bfm \psi}\in
L^2(\Omega)$. Then, by using \eqref{defrl} and \eqref{resolR}, we obtain
\begin{align*}
(R_\lambda{\bfm \varphi},{\bfm \psi})_2=&(PG_\lambda  P^*{\bfm
\varphi},{\bfm \psi})_2=(G_\lambda  P^*{\bfm \varphi},P^*{\bfm
\psi})=B_\lambda\big(G_\lambda  P^*{\bfm \psi},G_\lambda
P^*{\bfm \varphi}\big)\\=& B_\lambda\big(G_\lambda P^*{\bfm
\varphi},G_\lambda  P^*{\bfm \psi}\big)=( G_\lambda  P^*{\bfm \psi},P^*{\bfm
\varphi})=({\bfm
\varphi},R_\lambda{\bfm \psi})_2
\end{align*}
so that $R_\lambda$  is self-adjoint in $L^2(\Omega)$.

We now prove 2).  Consider ${\bfm \varphi}\in L^2(\Omega)$ satisfying  $ \lambda R_\lambda {\bfm \varphi}={\bfm \varphi}$ for some $\lambda>0$. We plug $ {\bfm g}= G_\lambda P^*{\bfm
\varphi}\in \mathds{W}^1$ into \eqref{resolR1}:
\begin{equation}\label{relationgl}
\lambda|R_\lambda{\bfm \varphi}|_2^2+\frac{1}{2}\int_{\Omega^+}{\bfm a}^+_{ij}\,\partial_i(G_\lambda
P^*{\bfm \varphi})\,\partial_j(G_\lambda P^*{\bfm
\varphi})\,d\mu^+=(PG_\lambda P^*{\bfm
\varphi},{\bfm \varphi})=(R_\lambda {\bfm \varphi},{\bfm \varphi})_2.
\end{equation}
Since $\lambda R_\lambda {\bfm \varphi}={\bfm \varphi}$, the right-hand side matches $(R_\lambda {\bfm \varphi},{\bfm \varphi})_2=\lambda|R_\lambda{\bfm \varphi}|_2^2$ so that the integral term in \eqref{relationgl} must vanish, that is $\int_{\Omega^+}{\bfm a}^+_{ij}\,\partial_i(G_\lambda
P^*{\bfm \varphi})\,\partial_j(G_\lambda P^*{\bfm
\varphi})\,d\mu^+=0$. From \eqref{unif_ellip}, we deduce $\partial
(G_\lambda P^*{\bfm \varphi})=0$. Thus, $G_\lambda P^*{\bfm \varphi}(0,\cdot)$ is  ${\cal
G}^* $-measurable. Moreover, we have $\lambda G_\lambda P^*{\bfm \varphi}(0,\cdot)=\lambda PG_\lambda P^*{\bfm \varphi}=\lambda R_\lambda {\bfm
\varphi}={\bfm
\varphi}$ so that ${\bfm \varphi}$ is ${\cal
G}^* $-measurable. Hence $ {\bfm
\varphi}=\M_1[{\bfm \varphi}]$. 

Conversely, we  assume ${\bfm \varphi}=\M_1[{\bfm \varphi}]$, which equivalently means that ${\bfm
\varphi}$ is $\mathcal{G}^*$ measurable. We define the
function ${\bfm u}:\Omega^+\to \R$ by ${\bfm u}(x_1,\omega)={\bfm \varphi}(\omega)$. It is obvious to check that ${\bfm u}$  belongs to $\mathds{W}^1$ and satisfies $\partial{\bfm u}=0 $. So $B_\lambda({\bfm u},\cdot)=(\cdot,\lambda P^*{\bfm \varphi})$ for any $\lambda>0$. This means ${\bfm u}=\lambda G_\lambda P^*{\bfm \varphi}$ in such a way that $\lambda R_\lambda{\bfm \varphi}=\lambda PG_\lambda P^*{\bfm \varphi}=P(\lambda G_\lambda P^*{\bfm \varphi})=P{\bfm u}={\bfm \varphi}$. %The same argument holds to prove: ${\bfm \varphi}=\M_1[{\bfm \varphi}]\Leftrightarrow{\bfm
%\varphi}\in {\rm Ker}(\lambda R_\lambda -{\rm I}).$

We prove 3). Consider ${\bfm \varphi}\in L^2(\Omega)$.  Since the relation \eqref{relationgl} is valid in great generality,  \eqref{relationgl} remains valid for such a function ${\bfm \varphi}$. Since the integral term in \eqref{relationgl} is nonnegative, we deduce $\lambda|R_\lambda{\bfm \varphi}|_2^2\leq (R_\lambda {\bfm \varphi},{\bfm \varphi})_2\leq |R_\lambda {\bfm \varphi}|_2| {\bfm \varphi}|_2$. Hence
$|\lambda R_\lambda{\bfm \varphi}|_2\leq |{\bfm \varphi}|_2$ for any $\lambda>0$. So the family $(\lambda
R_\lambda{\bfm \varphi})_\lambda $ is bounded in $L^2(\Omega)$ and we can extract
a subsequence, still indexed by $\lambda>0$, such that $(\lambda
R_\lambda{\bfm \varphi})_\lambda $ weakly converges in $L^2(\Omega) $
towards a function $\hat{\bfm \varphi}$.  Our purpose is now to establish that there is a unique possible weak limit $\hat{\bfm \varphi}=\M_1[{\bfm \varphi}] $ for the family $(\lambda
R_\lambda{\bfm \varphi})_\lambda $.

By multiplying the resolvent
relation $(\lambda-\mu)R_\lambda R_\mu{\bfm \varphi}=R_\mu{\bfm
\varphi}-R_\lambda{\bfm \varphi} $ by $\mu$ and passing to the limit as $ \mu\to 0$,
we get $\lambda R_\lambda \hat{\bfm \varphi}=\hat{\bfm \varphi}$. This latter relation implies (see above) that $\hat{\bfm \varphi}$ is ${\cal G}^*$-measurable. To prove $\hat{\bfm \varphi}=\M_1[{\bfm \varphi}] $, it just remains to establish the relation $({\bfm \varphi},{\bfm \psi})_2=(\hat{\bfm \varphi},{\bfm \psi})_2$ for every ${\cal G}^*$-measurable function ${\bfm \psi}\in L^2(\Omega)$. So we consider such a function ${\bfm \psi}$. Obviously, it satisfies the relations  $\M_1[{\bfm \psi}]={\bfm \psi}$ and $\lambda R_\lambda{\bfm \psi}={\bfm \psi}$ (see the above item 2). We deduce
$$({\bfm \varphi},{\bfm \psi})_2=({\bfm \varphi},\lambda R_\lambda {\bfm \psi})_2=\lim_{\lambda\to 0}(\lambda R_\lambda{\bfm \varphi},{\bfm
\psi})_2=(\hat{\bfm \varphi},{\bfm \psi})_2.$$ 
As a consequence, we have $\hat{\bfm \varphi}=\M_1[{\bfm \varphi}]$ and there is a unique possible limit for each weakly converging subsequence of the family $(\lambda R_\lambda{\bfm \varphi})_\lambda$. The whole family is therefore weakly converging in $L^2(\Omega)$. 

To establish the strong
convergence, it suffices to prove the convergence of the norms. As a weak limit,  $\hat{\bfm \varphi}$ satisfies the property $ |\hat{\bfm \varphi}|_2\leq
\liminf_{\lambda \to 0}|\lambda R_\lambda{\bfm \varphi}|_2$. Conversely, \eqref{relationgl} yields
 $$\limsup_{\lambda\to 0} | \lambda R_\lambda {\bfm
\varphi}|_2^2\leq \limsup_{\lambda\to 0} (\lambda R_\lambda {\bfm
\varphi},{\bfm \varphi})_2=(\hat{\bfm \varphi},{\bfm
\varphi})_2=|\hat{\bfm \varphi}|_2^2$$ and the
strong convergence follows.\qed

The remaining part of this section is concerned with the regularity properties of the operator $G_\lambda P^*$ (Propositions \ref{propregu} and \ref{max}) and may be omitted upon the first reading. Indeed, though they may appear a bit tedious, they are a direct adatation of existing results for the corresponding operators defined on $\bar{D}$ (not on $\Omega^+$). However, since we cannot quote proper references, we give the details.

Given ${\bfm u}\in L^2( \Omega^+)$, we shall say that ${\bfm
u}$ is a weakly differentiable if, for $i=1,\dots,d$, we can find some function
$\partial_i{\bfm u}\in L^2(\Omega^+)$ such that, for any
${\bfm g}\in \mathds{C}_c(\Omega^+)$:
$$\int_{\Omega^+}{\bfm u}\partial_i{\bfm g}\,d\mu^+= -\int_{\Omega^+}\partial_i{\bfm u}{\bfm g}\,d\mu^+.$$
It is straightforward to check that a function ${\bfm u}\in \mathds{W}^1$ is weakly differentiable.
For $k\geq 2$, the space $\mathds{W}^k$ is recursively defined as the set of
functions ${\bfm u}\in\mathds{W}^1 $ such that $\partial_i{\bfm u}$
is $k-1$ times weakly differentiable for $i=1,\dots,d$.

\begin{proposition}\label{wk}
If ${\bfm \varphi}$ belongs to $ {\cal C} $, then $G_\lambda P^*{\bfm
\varphi}\in \bigcap_{k=1}^\infty\mathds{W}^k$.
\end{proposition}

\noindent \textit{Proof of Proposition \ref{wk}.} The strategy is based on the well-known method
of difference quotients.  Our proof, adapted to the context of random media, is based on  \cite[Sect. 7.11 \& Th. 8.8]{gilbarg}. The properties of difference quotients in random media are summarized below (see e.g. \cite[Sect. 5]{rhodes:06}):

i) for $j=2,\dots,d$, $r\in\R\setminus\{0\}$ and ${\bfm g}\in \mathds{C}_c(\Omega_+)$, we define $$\Delta_r^j{\bfm g}(x_1,\omega)=\frac{1}{r}({\bfm g}(x_1,\tau_{re_j}\omega)-{\bfm g}(x_1,\omega)).$$

ii) for each $r\in\R\setminus\{0\}$ and ${\bfm g}\in \mathds{C}_c(\Omega_+)$, we define $$\Delta_r^1{\bfm g}=\frac{1}{r}({\bfm g}(x_1+r,\omega)-{\bfm g}(x_1,\omega)).$$

iii) for any $j=1,\dots,d$, $r\in\R\setminus\{0\}$ and ${\bfm g},{\bfm h}\in \mathds{C}_c(\Omega_+)$, the discrete integration by parts holds
$$\int_{\Omega^+}\Delta_r^j{\bfm g}{\bfm h}\,d\mu^+=- \int_{\Omega^+}{\bfm g}\Delta_{-r}^j{\bfm h}\,d\mu^+$$ provided that $r$ is small enough to ensure that $\Delta_r^j{\bfm g}$ and $\Delta_r^j{\bfm h}$ belong to $ \mathds{C}_c(\Omega_+)$.

iv) for any $j=1,\dots,d$, $r\in\R\setminus\{0\}$ and ${\bfm g}\in \mathds{C}_c(\Omega_+)$ such that $\Delta_r^j{\bfm g}\in\mathds{C}_c(\Omega_+)$, we have
$$\int_{\Omega^+}|\Delta_r^j{\bfm g}|^2\,d\mu^+ \leq \int_{\Omega^+}|\partial_j{\bfm g}|^2\,d\mu^+.$$

Up to the end of the proof, the function $ G_\lambda P^*{\bfm \varphi}$ is denoted by ${\bfm u}$. The strategy consists in differentiating the resolvent equation $B_\lambda({\bfm u},\cdot)=(\cdot,P^*{\bfm \varphi}) $ to prove that the derivatives of ${\bfm u} $ equations of the same type. For $p=2,\dots,d$, it raises no difficulty to adapt the method explained in  \cite[Sect. 5]{rhodes:06} and prove that the "tangential derivatives" $\partial_p {\bfm u}$ 
belongs to $\mathds{W}^1 $ and solves the equation
\begin{equation}\label{equ:der}
B_\lambda(\partial_p{\bfm u},\cdot)=(\cdot,P^*D_p{\bfm
\varphi})-F_p(\cdot),
\end{equation}
 where $F_p:\mathds{W}^1\to\R $ is
defined by $$F_p({\bfm g})=(1/2)\int_{\Omega^+}D_p{\bfm
a}^+%+D_p{\bfm H}^+)
_{ij}\,\partial_i{\bfm u}\partial_j{\bfm
g}\,d\mu^+.$$
In particular, $\partial_{ij}{\bfm u}\in L^2(\Omega^+;\mu^+)$ for $(i,j)\not = (1,1)$. We let the reader check the details.

The main difficulty lies in  the "normal derivative" $\partial_{1}{\bfm u} $: we have to prove that  $\partial_{1}{\bfm u} $ is weakly differentiable. Actually, it just remains to prove that there exists a function $\partial^2_{11}{\bfm u}\in L^2(\Omega^+;\mu^+) $ such that $\forall{\bfm
g}\in\mathds{C}_c(\Omega_+)$:
\begin{equation}\label{d1u}
\int_{\Omega^+}\partial^2_{11}{\bfm u}{\bfm g}\,d\mu^+ =- \int_{\Omega^+}\partial_{1}{\bfm u}\partial_1{\bfm g}\,d\mu^+ .
\end{equation}
To that purpose, we plug a generic function ${\bfm
g}\in\mathds{C}_c(\Omega_+)$ into the resolvent equation \eqref{resolR1}. The boundary terms $(P^*{\bfm \varphi},{\bfm g})=({\bfm \varphi},{\bfm g}(0,\cdot))_2$ and $\lambda (P{\bfm u},P{\bfm g})_2 $ vanish and we obtain:
$$ \sum_{i,j=1}^d\int_{\Omega^+}{\bfm a}^+
_{ij}\partial_i{\bfm
u}\partial_j{\bfm
g}\,d\mu^+=0.$$
We isolate the term corresponding to $i=1$ and $j=1$ to obtain  (remind that ${\bfm a}_{11}=1$)
% and $\bar{{\bfm H}}_{11}=0$)
\begin{align*}
\int_{\Omega^+}\partial_1{\bfm u}\partial_1{\bfm
g}\,d\mu^+=&-\sum_{(i,j)\not=(1,1)}\int_{\Omega^+}{\bfm a}^+%+\bar{{\bfm H}}^+)
_{ij}\partial_i{\bfm
u}\partial_j{\bfm
g}\,d\mu^+\\=&\sum_{(i,j)\not=(1,1)}\int_{\Omega^+}\partial_j{\bfm a}^+%+\bar{{\bfm H}}^+)
_{ij}\partial_i{\bfm
u}{\bfm
g}\,d\mu^++\sum_{(i,j)\not=(1,1)}\int_{\Omega^+}{\bfm a}^+%+\bar{{\bfm H}})
_{ij}\partial^2_{ij}{\bfm
u}{\bfm g}\,d\mu^+.
\end{align*}
Since $\partial_{ij}{\bfm u}\in L^2(\Omega^+;\mu^+)$ for $(i,j)\not = (1,1)$, we deduce that $$\int_{\Omega^+}\partial_1{\bfm u}\partial_1{\bfm
g}\,d\mu^+\leq
C\big(\int_{\Omega^+}{\bfm g}^2\,d\mu^+\big)^{1/2}$$ for some positive constant $C$. So the mapping ${\bfm g}\in\mathds{C}_c(\Omega_+)\mapsto \int_{\Omega^+}\partial_1{\bfm u}\partial_1{\bfm
g}\,d\mu^+$ is $L^2(\Omega^+;\mu^+) $-continuous and there exists a unique function denoted by $\partial^2_{11}{\bfm u}$ such that \eqref{d1u} holds. As a consequence, $\partial_1{\bfm u}$ is weakly differentiable, that is  ${\bfm u}\in \mathds{W}^2$. Note that \eqref{equ:der} only involves the
functions ${\bfm a},{\bfm \varphi}$ and their
derivatives in such a way that we can iterate the argument in
differentiating \eqref{equ:der} and so on. So it is clear that the proof can be completed recursively.\qed

\vspace{3mm}
\noindent \textit{Proof of Proposition \ref{propregu}.} The function ${\bfm u}$ still stands for $G_\lambda P^*{\bfm \varphi}$. From Proposition \ref{wk}, we have ${\bfm u}\in \bigcap_{k=1}^\infty\mathds{W}^k$ and it is plain to deduce that $\mu$ a.s. the trajectories of ${\bfm u}$ are smooth and 
\begin{equation}\label{diffu}
\forall x=(x_1,y)\in \bar{D},\quad \partial_{x_i}\tilde{u}_\omega(x)=\partial_i{\bfm u}(x_1,\tau_{(0,y)}\omega) .
\end{equation}
 We let the reader check that point (it is a straightforward adaptation of the fact that an infinitely weakly differentiable function $f:\bar{D}\to \R$ is smooth).

It remains to prove that $\tilde{u}_\omega$ solves
\eqref{prob:halfphi}. To begin with, we state the following lemma
\begin{lemma}\label{jonc}
For each function ${\bfm v}\in \mathds{W}^1$, we define $\tilde{v}_\omega:(x_1,y)\in\bar{D}\mapsto {\bfm v}(x_1,\tau_{(0,y)}\omega) $. Then for every $\varrho\in
C_c^\infty(\bar{D})$ and ${\bfm \psi}\in {\cal C}$ we have:
\begin{align*}
\M\big[{\bfm \psi}(\omega)\int_{\partial D}(\lambda \tilde{v}_\omega(y)
 -{\bfm \gamma}_i(\tau_y\omega)\partial_{x_i}\tilde{v}_\omega(y))\varrho(y)\,dy\big]= B_\lambda({\bfm v},{\bfm
\psi}\ast \varrho)+ \M\big[{\bfm \psi}(\omega)\int_{\bar{D}} L^\omega \tilde{v}_\omega
(x)\varrho(x)\,dx\big]
\end{align*}
 where the function ${\bfm
\psi}\ast \varrho:\Omega^+\to \R$ belongs to $\mathds{W}^1$ and is defined by:
$${\bfm
\psi}\ast \varrho (x_1,\omega)=\int_{\R^{d-1}}\varrho(x_1,-y){\bfm \psi}(\tau_y\omega)\,dy .$$ 
\end{lemma}

 Let us consider
$\varrho\in C_c^\infty(\bar{D})$, ${\bfm \psi}\in {\cal
C}$. We first point out that 
\begin{align*}
B_\lambda({\bfm u},{\bfm
\psi}\ast \varrho)=&({\bfm
\psi}\ast \varrho,P^*{\bfm \varphi}) =\M\big[{\bfm \varphi}(\omega)\int_{\R^{d-1}}\!\!\!{\bfm \psi}(\tau_{(0,y)}\omega)\varrho(0,-y)\,dy\big]=
\M\big[{\bfm \psi}(\omega)\int_{\partial
D}\!\!\!{\bfm \varphi}(\tau_y\omega)\varrho(y)\,dy\big] .
\end{align*}
Then, by using  Lemma \ref{jonc} and the above relation, we obtain
\begin{align*}
\M\big[{\bfm \psi}(\omega)\int_{\partial D} \big(\lambda
\tilde{u}_\omega(y)
-{\bfm \gamma}_i(\tau_y\omega)\partial_{x_i}\tilde{u}_\omega(y)-{\bfm \varphi}(\tau_y\omega)\big)\varrho(y)\,dy\big]=\M\big[{\bfm \psi}(\omega)\int_{\bar{D}} L^\omega \tilde{u}_\omega
(x)\varrho(x)\,dx\Big].
\end{align*}
%= & \lambda(P{\bfm u},P({\bfm
%\chi}\ast \check{\varrho}))_2+(1/2)\int_{\Omega^+}{\bfm a}^+%+\bar{{\bfm H}}^+)_{ij}\partial_i
%{\bfm u}
%\partial_j({\bfm
%\chi}\ast \check{\varrho})\,d\mu^++ \M\big[{\bfm \chi}(\omega)\int_{\bar{D}} L^\omega \tilde{u}_\omega
%(x)\varrho(x)\,dx\big]\\
%= & B_\lambda({\bfm u},{\bfm
%\chi}\ast \check{\varrho})+ \M\big[{\bfm \chi}(\omega)\int_{\bar{D}} L^\omega \tilde{u}_\omega
%(x)\varrho(x)\,dx\big]=({\bfm
%\chi}\ast \check{\varrho},P^*{\bfm \varphi})+ \M\big[{\bfm \chi}(\omega)\int_{\bar{D}} L^\omega \tilde{u}_\omega
%(x)\varrho(x)\,dx\big]\\=&
%\int_{\bar{D}} L^\omega \tilde{u}_\omega
%(x)\varrho(x)\,dx\Big)\big].
Since the above relation is valid for any ${\bfm \psi}\in {\cal C}$, we deduce that $\mu$ a.s. we have
$$\int_{\partial D} \big(\lambda
\tilde{u}_\omega(y)
-{\bfm \gamma}_i(\tau_y\omega)\partial_{x_i}\tilde{u}_\omega(y)-{\bfm \varphi}(\tau_y\omega)\big)\varrho(y)\,dy= \int_{\bar{D}} L^\omega \tilde{u}_\omega
(x)\varrho(x)\,dx.$$ 
By choosing in turn a generic function $\varrho $ vanishing or not on the boundary, we deduce that $\mu$ a.s. we have: $L^\omega \tilde{u}_\omega=0$ on $D$ and $\lambda \tilde{u}_\omega(y)
-{\bfm \gamma}_i(\tau_y\omega)\partial_{x_i}\tilde{u}_\omega(y)={\bfm \varphi}(\tau_y\omega)$ for 
$y\in \partial D$.\qed

\vspace{2mm}
\noindent {\it Proof of Lemma \ref{jonc}.} %Since the proof of the above lemma relies does not raise any difficulty, we just give a hint of proof. 
First apply the Green formula \eqref{green2} (with $V=0$ and $\epsilon=1 $):
\begin{align*}
\int_{\partial D}(\lambda \tilde{v}_\omega(y)
 -{\bfm \gamma}_i(\tau_y\omega)\partial_{x_i}\tilde{v}_\omega(y))\varrho(y)\,dy =& \int_{\partial D}\lambda \tilde{v}_\omega(y)\varrho(y)\,dy +\frac{1}{2}\int_{\bar{D}}{\bfm a}_{ij}(\tau_x\omega)\partial_{x_i} \tilde{v}_\omega(x)\partial_{x_j}\rho(x)\,dx\\&+\int_{\bar{D}} L^\omega \tilde{v}_\omega 
(x)\varrho(x)\,dx.
\end{align*}
Then we multiply the above relation by $ {\bfm \psi}$ and integrate with respect to $\M$. By using the invariance of $\mu$ under translations, we have
\begin{align*}
\M\big[ {\bfm \psi}(\omega)\int_{\partial D}\lambda \tilde{v}_\omega(y)\varrho(y)\,dy\big]= &\lambda\M\big[\int_{\partial D} {\bfm \psi}(\omega) {\bfm v}(0,\tau_y\omega)\varrho(y)\,dy\big]= \lambda\M\big[{\bfm v}(0,\omega)\int_{\partial D} {\bfm \psi}(\tau_{-y}\omega) \varrho(y)\,dy\big] \\
=&\lambda (P{\bfm v},P{\bfm
\psi}\ast \varrho)_2.
\end{align*}
With similar arguments and \eqref{diffu}, we prove
$$\frac{1}{2}\int_{\bar{D}}{\bfm a}_{ij}(\tau_x\omega)\partial_{x_i} \tilde{v}_\omega(x)\partial_{x_j}\rho(x)\,dx=\M\int_{\R_+}{\bfm a}_{ij}^+\partial_i{\bfm v}\partial_j {\bfm
\psi}\ast \varrho\,d\mu^+.$$
The lemma follows.\qed

\vspace{3mm}
\noindent \textit{Proof of Proposition \ref{max}.} We adapt the Stampacchia truncation method. More precisely, we introduce a function $H:\R\to \R$ of class $C^1(\R)$ such that 
$$i)\forall s \in \R,\,\,|H'(s)|\leq C,\quad ii)\forall s >0,\,\, H'(s)>0,\quad  iii)\forall s \leq 0,\,\, H'(s)=0.$$
We define $K=|{\bfm \varphi}|_\infty/\lambda$ and ${\bfm u}_\lambda=G_\lambda P^*{\bfm \varphi}$. We let the reader check that $H({\bfm u}_\lambda-K)\in\mathds{W}^1$. Then we plug ${\bfm g}=H({\bfm u}_\lambda-K) $ into \eqref{resolR1} and we obtain:
$$ \lambda(P{\bfm u}_\lambda,PH({\bfm u}_\lambda-K))_2+\frac{1}{2}\int_{\Omega^+}{\bfm a}^+_{ij}\,\partial_i{\bfm u}_\lambda\,\partial_j
{\bfm u}_\lambda H'({\bfm u}_\lambda-K) \,d\mu^+=(PH({\bfm u}_\lambda-K),{\bfm \varphi})_2.$$ 
By subtracting the term $\lambda (K,H(P{\bfm u}_\lambda-K))_2 $ in each side of the above equality, we obtain:
$$ \lambda(P{\bfm u}_\lambda-K,H(P{\bfm u}_\lambda-K))_2+\frac{1}{2}\int_{\Omega^+}{\bfm a}^+_{ij}\,\partial_i{\bfm u}_\lambda\,\partial_j
{\bfm u}_\lambda H'({\bfm u}_\lambda-K) \,d\mu^+=(H(P{\bfm u}_\lambda-K),{\bfm \varphi}-\lambda K)_2.$$
Observe that the right-hand side is negative since ${\bfm \varphi}-\lambda K\leq 0$ and $H(s)\geq 0$ for any $s\geq 0$. Furthermore, the left-hand side is positive since $H'(s)\geq 0$   and $sH(s)\geq 0$ for $s\in \R$. We deduce that both terms of the left-hand side reduce to $0$. The relation $ \lambda(P{\bfm u}_\lambda-K,H(P{\bfm u}_\lambda-K))_2=0$ and the properties of $H$ ($sH(s)\geq 0$ for $s\in \R$ and $sH(s)> 0$ for $s>0$) ensure that $P{\bfm u}_\lambda -K\leq 0$, that is
\begin{equation}\label{stamp1}
PH({\bfm u}_\lambda-K)=0.
\end{equation}
The relation $\frac{1}{2}\int_{\Omega^+}{\bfm a}^+_{ij}\,\partial_i{\bfm u}_\lambda\,\partial_j
{\bfm u}_\lambda H'({\bfm u}_\lambda-K) \,d\mu^+=0$ and \eqref{unif_ellip} prove that $|\partial {\bfm u}_\lambda|^2\,H'({\bfm u}_\lambda-K)=0$  $\mu^+$ a.s.. In particular, we deduce that 
\begin{equation}\label{stamp2}
\partial \big( H({\bfm u}_\lambda-K)\big)=0 .
\end{equation}
By gathering \eqref{stamp1} and \eqref{stamp2}, we deduce $N(H({\bfm u}_\lambda-K))=0$ (recall the definition of $N$ in \eqref{norm}). So $H({\bfm u}_\lambda-K)=0$ and this means ${\bfm u}_\lambda\leq K$.\qed

%%%%%%%%%%%%%%%%%%%%%%%%%%%%%%%%%%%%%%%%%%%%%%%%%%%%%%%%%%%%%%%%%%%%%%%%%%%%%%%%%%%%%%%%%%%%%%%%%%%%%%
\section{Proofs of subsection \ref{ergth}}
%%%%%
%%%%

\noindent \textit{Proof of Theorem \ref{theorem_ergodic}.} We first suppose that ${\bfm f}$ belongs to $\mathcal{C}$. Even
if it means replacing ${\bfm f}$ by ${\bfm  f}-\M[{\bfm
f}]$, it is enough to treat the case $\M[{\bfm f}]=0$. We consider the solution ${\bfm
v}_\lambda\in L^2(\Omega)\cap {\rm Dom}({\bfm L})$ to the resolvent equation
\begin{equation}\label{eq_resolvent}
\lambda {\bfm v}_\lambda-{\bfm L}{\bfm v}_\lambda={\bfm f}.
\end{equation}
For the same reason as in the proof of Proposition \ref{propregu}, $\mu$ a.s. the function $\vartheta:x\in\R^d\mapsto
{\bfm v}_\lambda(\tau_x\omega) $ satisfies   $\lambda
\vartheta(x)-L^\omega\vartheta (x)={\bfm f}(\tau_x\omega)  $
$x\in\R^d$. So $\vartheta$ is smooth \cite[Th. 6.17]{gilbarg}.
Applying the It\^o formula to the function $x\mapsto
{\bfm v}_\lambda(\tau_x\omega) $ then yields
\begin{equation*}
 \begin{split}
d{\bfm v}_\lambda (\tau_{X^\varepsilon_t/\varepsilon}\omega)
=& \varepsilon^{-1}D_i
{\bfm v}_\lambda{\bfm \sigma}_{ij}(\tau_{X^\varepsilon_t/\varepsilon}\omega)dB^{*j}_t -\varepsilon^{-1}\partial_{x_i}V(X^\varepsilon_t){\bfm a}_{ij}
 D_j{\bfm v}_\lambda(\tau_{X^\varepsilon_t/\varepsilon}\omega)\,dt\\&+\varepsilon^{-2}
 {\bfm L} {\bfm v}_\lambda (\tau_{X^\varepsilon_t/\varepsilon}\omega)\,dt+\varepsilon^{-1}
 D_i{\bfm v}_\lambda {\bfm \gamma}_i(\tau_{X^\varepsilon_t/\varepsilon}\omega)\,dK^\varepsilon_t.
 \end{split}
\end{equation*} In the above expression, we  replace ${\bfm L}{\bfm v}_\lambda$ by $\lambda {\bfm v}_\lambda-{\bfm f} $, multiply both sides of the equality by $\varepsilon^2$ and isolate the term  ${\bfm f}(\tau_{X^\varepsilon_t/\varepsilon}\omega)\,dt$. We obtain
\begin{align}
 \int_0^t \! {\bfm f}(\tau_{X^\varepsilon_r/\varepsilon}\omega\!)dr\!=&
\varepsilon\!\int_0^t\! D_i
{\bfm v}_\lambda{\bfm \sigma}_{ij}(\tau_{X^\varepsilon_r/\varepsilon}\omega)dB^{*j}_r\!
 -\!\varepsilon^2\!\big( {\bfm v}_\lambda(\tau_{X^\varepsilon_t/\varepsilon}\omega)\!-\!
{\bfm v}_\lambda(\tau_{X^\varepsilon_0/\varepsilon}\omega)\big)\!+\!\int_0^t\!\lambda
{\bfm v}_\lambda (\tau_{X^\varepsilon_r/\varepsilon}\omega)dr\nonumber\\&\!+\!\varepsilon\!\int_0^t\!
 D_i{\bfm v}_\lambda {\bfm \gamma}_i(\tau_{X^\varepsilon_t/\varepsilon}\omega)dK^\varepsilon_r\!-\!\varepsilon\!\int_0^t\!\partial_{x_i}V(X^\varepsilon_r){\bfm a}_{ij}D_j
{\bfm v}_\lambda (\tau_{X^\varepsilon_r/\varepsilon}\omega)dr\label{sum}\\
 \stackrel{\mathrm{def}}{=}&\Delta^{1,\varepsilon,\lambda}_t
-\Delta^{2,\varepsilon,\lambda}_t +\Delta^{3,\varepsilon,\lambda}_t+
\Delta^{4,\varepsilon,\lambda}_t-
\Delta^{5,\varepsilon,\lambda}_t.\nonumber
\end{align}
Let us investigate the quantities $\Delta^{1,\varepsilon,\lambda}$, $\Delta^{2,\varepsilon,\lambda}$, $\Delta^{3,\varepsilon,\lambda}$, $\Delta^{4,\varepsilon,\lambda}$ and $\Delta^{5,\varepsilon,\lambda}$. Using the Doob inequality and Lemma \ref{lem_invariant}, we have:
$$\bar{\E}^{\varepsilon *}\big[\sup_{0\leq t \leq T}|\Delta^{1,\varepsilon,\lambda}_t|^2\big]\leq 4\varepsilon^2T\M_D^*\big[|D_i
{\bfm v}_\lambda{\bfm \sigma}_{ij}|^2\big]\leq C\varepsilon^2|D
{\bfm v}_\lambda|_2^2$$ for some positive constant $C$ only depending on $T$ and $|{\bfm \sigma}|_\infty $. Hence $ \bar{\E}^{\varepsilon *}\big[\sup_{0\leq t \leq
T}|\Delta^{1,\varepsilon,\lambda}_t|^2\big]\to 0$ as $\varepsilon\to
0$, for each fixed $ \lambda>0$. Similarly, by using the boundedness of ${\bfm a},{\bfm \gamma},\partial_xV$, we can prove 
$$\bar{\E}^{\varepsilon *}\big[\sup_{0\leq t \leq
T}|\Delta^{4,\varepsilon,\lambda}_t|+\sup_{0\leq t \leq
T}|\Delta^{5,\varepsilon,\lambda}_t|^2\big]\to 0,\quad \text{as }\varepsilon\to 0.$$
From Lemma \ref{maxU}, ${\bfm v}_\lambda$ is bounded by $|{\bfm f}|_\infty/\lambda $. We deduce
$$ \bar{\E}^{\varepsilon *}\big[\sup_{0\leq t \leq T}|\Delta^{2,\varepsilon,\lambda}_t|^2\big]\leq 4\varepsilon^4|{\bfm f}|_\infty^2 \lambda^{-2}\to 0,\quad \text{as }\varepsilon\to 0.$$ 
By taking the $\limsup_{\epsilon\to 0}$ in \eqref{sum} and by using the convergences of $\Delta^{1,\varepsilon,\lambda},\Delta^{2,\varepsilon,\lambda},\Delta^{4,\varepsilon,\lambda},\Delta^{5,\varepsilon,\lambda} $ towards $0$, we deduce
\begin{align*}
\limsup_{\varepsilon\to 0}\bar{\E}^{\varepsilon *}\big[\sup_{0\leq t
\leq T}|\int_0^t
{\bfm f}(\tau_{X^\varepsilon_r/\varepsilon}\omega)\,dr|\big]&\leq\limsup_{\varepsilon\to
0}\bar{\E}^{\varepsilon *}\big[\sup_{0\leq t \leq
T}|\Delta^{3,\varepsilon,\lambda}_t|\big].
\end{align*}
Furthermore, from  Lemma \ref{lem_invariant}, we have
$$\limsup_{\varepsilon\to
0}\bar{\E}^{\varepsilon *}\big[\sup_{0\leq t \leq
T}|\Delta^{3,\varepsilon,\lambda}_t|\big]\leq \limsup_{\varepsilon\to
0}\int_0^T\bar{\E}^{\varepsilon *}\big[|\lambda
{\bfm v}_\lambda (\tau_{X^\varepsilon_r/\varepsilon}\omega)|\big]\,dr= T|\lambda{\bfm
v}_\lambda|_1\leq T|\lambda{\bfm
v}_\lambda|_2.
 $$
From Proposition \ref{ergprob}, we have $|\lambda{\bfm
v}_\lambda|_2\rightarrow 0$ as $\lambda$ goes to $0$. So it just remains to choose $\lambda$ small enough to complete the
proof in the case of a smooth function ${\bfm f}\in\mathcal{C}$. The general
case follows from the density of $ \mathcal{C}$ in
$L^1(\Omega)$ and Lemma \ref{lem_invariant}.\qed

\vspace{2mm}

\noindent \textit{Proof of Theorem \ref{theorem_erglocal}.} Once again, from Lemma \ref{lem_invariant} and density
arguments, it is sufficient to consider the case of a smooth
function ${\bfm f}\in {\cal C}$. Even if it means replacing ${\bfm f}$ with ${\bfm f}-\M_1[{\bfm f}]$, it is enough to consider the case $\M_1[{\bfm f}]=0$. Let us define, for any $\lambda>0$,
${\bfm u}_\lambda= G_\lambda P^*{\bfm f}$ and ${\bfm
f}_\lambda=R_\lambda {\bfm f}$, the definitions of which are given in Section \ref{ergo} (boundary ergodic problems). We still use the notation
$\tilde{u}_\omega^\lambda(x)={\bfm
u}_\lambda(x_1,\tau_{(0,y)}\omega) $ for any
$x=(x_1,y)\in \bar{D}$. We remind the reader that the main regularity properties of the function $ \tilde{u}_\omega^\lambda$ are summarized in Proposition \ref{propregu}. In particular, $\mu$ a.s., the mapping $x\mapsto \tilde{u}_\omega^\lambda(x)$ is smooth 
and we can apply the It\^o formula: 
\begin{align}
d\big(\varepsilon\tilde{u}_\omega^\lambda(X^\varepsilon_t/&\varepsilon)\big)=\big[\varepsilon^{-1}L^\omega\tilde{u}_\omega^\lambda(X^\varepsilon_t/\varepsilon)-\partial_{x_j}V(X^\varepsilon_t){\bfm a}_{ij}(\tau_{X^\varepsilon_t/\varepsilon}\omega)\partial_{x_i}\tilde{u}_\omega^\lambda(X^\varepsilon_t/\varepsilon)\big]\,dt\nonumber\\&\quad +\partial_{x_i}
\tilde{u}_\omega^\lambda(X^\varepsilon_t/\varepsilon){\bfm \sigma}_{ij}(\tau_{X^\varepsilon_t/\varepsilon}\omega)\,dB_t^{*j}+{\bfm \gamma}_i(\tau_{X^\varepsilon_t/\varepsilon}\omega)\partial_{x_i}\tilde{u}_\omega^\lambda(X^\varepsilon_t/\varepsilon)\,dK^\varepsilon_t\label{itou}
\end{align}
In the above expression, we use the relation $ L^\omega\tilde{u}_\omega^\lambda=0$ inside $D$. Furthermore, since ${\bfm \gamma}_i\partial_{x_i} {u}_\omega^\lambda(x)=\lambda {\bfm f}_\lambda(\tau_x\omega)-{\bfm f}(\tau_x\omega)$ on $\partial D$ and $dK^\varepsilon_t=\one_{\partial D}(X^\varepsilon_t)dK^\varepsilon_t$, we deduce $${\bfm \gamma}_i(\tau_{X^\varepsilon_t/\varepsilon}\omega)\partial_{x_i}\tilde{u}_\omega^\lambda(X^\varepsilon_t/\varepsilon)\,dK^\varepsilon_t=(\lambda {\bfm f}_\lambda-{\bfm f})(\tau_{X^\varepsilon_t/\varepsilon}\omega)\,dK^\varepsilon_t .$$ Hence, \eqref{itou} yields
\begin{align}
\int_0^t{\bfm f}(\tau_{X^\varepsilon_r/\varepsilon}\omega)\,dK^\varepsilon_r=&-\Big(\varepsilon\tilde{u}_\omega^\lambda(X^\varepsilon_t/\varepsilon)-\varepsilon\tilde{u}_\omega^\lambda(X^\varepsilon_0/\varepsilon)\Big)-\int_0^t\partial_{x_j}V(X^\varepsilon_r){\bfm a}_{ij}(\tau_{X^\varepsilon_r/\varepsilon}\omega)\partial_{x_i}\tilde{u}_\omega^\lambda(X^\varepsilon_r/\varepsilon)\,dr\nonumber\\&+\int_0^t\partial_{x_i}
\tilde{u}_\omega^\lambda(X^\varepsilon_r/\varepsilon){\bfm \sigma}_{ij}(\tau_{X^\varepsilon_r/\varepsilon}\omega)\,dB_r^{*j}+\int_0^t\lambda
{\bfm f}_\lambda(\tau_{X^\varepsilon_r/\varepsilon}\omega)\,dK^\varepsilon_r\nonumber\\
\label{eq_itoergtheo}
\equiv & -\Delta^{1,\varepsilon}_t-
\Delta^{2,\varepsilon}_t+\Delta^{3,\varepsilon}_t+\Delta^{4,\varepsilon}_t.
\end{align}
The next step of the proof is to prove that $\Delta^{1,\varepsilon},\Delta^{2,\varepsilon},\Delta^{3,\varepsilon}$ converge to $0$ as $\varepsilon$ goes to $0$ for each fixed $\lambda>0 $. Clearly, from Proposition \ref{max}, we have $$\bar{\E}^{\varepsilon *}\big[\sup_{0\leq t \leq T}|\Delta^{1,\varepsilon}_t|^2\big]\leq 4\varepsilon^2|{\bfm u}_\lambda|^2_{L^\infty(\Omega^+)}\xrightarrow[\varepsilon\to 0]{} 0.$$
%From the identity ${\bfm u}_\lambda(x_1,\omega)-{\bfm
%u}_\lambda(0,\omega)=\int_0^{x_1}\partial_1{\bfm
%u}_\lambda(r,\omega)\,dr$ and the Cauchy-Schwarz inequality, we obtain
%$|{\bfm u}_\lambda(x_1,\cdot)|_2^2\leq 2|{\bfm
%u}_\lambda(0,\cdot)|_2^2+2x_1\int_0^{x_1}|\partial_1{\bfm
%u}_\lambda(r,\cdot)|_2^2\,dr$ so that
%\begin{align}\label{estint}
%  \M_{D}[|{\bfm
%  u}_\lambda(\cdot/\varepsilon,\cdot)|^2] \leq &2|{\bfm
%u}_\lambda(0,\cdot)|_2^2+\frac{2\M_{D}[x_1]}{\varepsilon}\int_{\Omega^+}|\partial_1{\bfm
%  u}_\lambda|^2\,d\mu^+\nonumber\\ \leq &  2|P{\bfm
%u}_\lambda|_2^2+\frac{2M}{\varepsilon}\int_{\Omega^+}|\partial_1{\bfm
%  u}_\lambda|^2\,d\mu^+
%\end{align}
Let us now focus on $\Delta^{2,\varepsilon}_t$. We use the
boundedness of $\partial_{x_j} V,{\bfm a}_{ij}$ ($1 \leq i,j\leq d$) and
Lemma \ref{lem_invariant}: $$\bar{\E}^{\varepsilon *}\big[\sup_{0\leq t \leq T}|\Delta^{2,\varepsilon}_t|^2\big]\leq 
 T|{\bfm a}|^2_\infty \times\sup_{\bar{D}}|\partial_xV| ^2 \times\M_D^*\big[|\partial_x \tilde{u}^\lambda_\omega(\cdot/\varepsilon)|^2 \big].$$  Furthermore
\begin{align}
 \M_D^*\big[|\partial_x \tilde{u}^\lambda_\omega(\cdot/\varepsilon)|^2 \big]\!=&\M\int_{(x_1,y)\in\bar{D}}|\partial {\bfm u}_\lambda(x_1/\epsilon,\tau_{y/\varepsilon}\omega)|^2e^{-2V(x_1,y)}\,dx_1dy\nonumber\\
 =&\M\!\int_{\R_+}\!|\partial {\bfm
u}_\lambda(x_1/\varepsilon,\omega)|^2\Big(\int_{\R^{d-1}}e^{-2V(x_1,y)}\,dy\Big)dx_1.\label{nimporte}
\end{align}
We point out that the function $V$ given by \eqref{defV} satisfies
\begin{equation}\label{V3}
%\int_{\bar{D}}x_1e^{-2V(x)}\,dx=M<+\infty\,\,\text{ and }\,\,
S\stackrel{def}{=}\sup_{x_1\geq 0}\int_{\R^{d-1}}e^{-2V(x_1,y)}\,dy<+\infty.
\end{equation}
By gathering \eqref{V3} and \eqref{nimporte} and by making the change of variables $u=x_1/\varepsilon$, we deduce that $ \M_D^*\big[|\partial_x \tilde{u}^\lambda_\omega(\cdot/\varepsilon)|^2 \big]$ is not greater than $\varepsilon S
\!\int_{\Omega^+}\!|\partial {\bfm
u}_\lambda|^2d\mu^+$. So $\bar{\E}^{\varepsilon *}\big[\sup_{0\leq t \leq T}|\Delta^{2,\varepsilon}_t|^2\big]$  converges to $0$ as $\epsilon\to 0$.
By combining the same argument with the Doob inequality, we prove that $\bar{\E}^{\varepsilon *}\big[\sup_{0\leq t \leq T}| \Delta^{3,\varepsilon}_t|^2\big]\to 0$ as
$\varepsilon\to 0$. 

So, taking the $\limsup_{\epsilon\to 0} $ in \eqref{eq_itoergtheo} and using the above convergences yields
$$\limsup_{\varepsilon\to 0}\bar{\E}^{\varepsilon *}\big[\sup_{0\leq t \leq T}|\int_0^t{\bfm f}(\tau_{X^\varepsilon_r/\varepsilon}\omega)\,dK^\varepsilon_r|\big]\leq \limsup_{\varepsilon\to 0}\bar{\E}^{\varepsilon *}\big[\int_0^T|\lambda{\bfm f}_\lambda(\tau_{X^\varepsilon_r/\varepsilon}\omega)|\,dK^\varepsilon_r\big].$$
By using Lemma \ref{lem_invariant} in the right-hand side of the previous inequality, we deduce, for any $\lambda>0$,
$$\limsup_{\varepsilon\to 0}\bar{\E}^{\varepsilon *}\big[\sup_{0\leq t \leq T}|\int_0^t{\bfm f}(\tau_{X^\varepsilon_r/\varepsilon}\omega)\,dK^\varepsilon_r|\big]\leq T\M^*_{\partial D}[|\lambda{\bfm f}_\lambda|]
=
 |\lambda{\bfm f}_\lambda|_1T\int_{\partial D}e^{-2V(x)}\,dx\leq S T |
\lambda{\bfm f}_\lambda|_2.$$
From Proposition \ref{erglocal} item 3, we can choose $\lambda$ small enough so
as to make the latter term arbitrarily small. So we complete the proof.\qed

\noindent \textit{Proof of Theorem \ref{main_erg}.} 1) From \eqref{girsanov}, we only have to check that \eqref{main_erg1} holds under $\bar{\P}^{\varepsilon *}$. This follows from Theorem \ref{theorem_ergodic} and the estimate (obtained with Lemma \ref{lem_invariant}) $$\lim_{\varepsilon\to 0}\bar{\E}^{\varepsilon *}\Big[\sup_{0 \leq t \leq
T}|\int_0^t({\bfm f}_\varepsilon-{\bfm f})(\tau_{X^{\varepsilon}_r/\varepsilon}\omega)\,dr|\Big]\leq T|{\bfm f}_\varepsilon-{\bfm f}|_1.$$
The same argument holds for \eqref{main_erg2}.\qed

%%%%%%%%%%%%%%%%%%%%%%%%%%%%%%%%%%%%%%%%%%%%%%%%%
\section{Proofs of subsection \ref{corr}}\label{app:corr}
%%%%%%%%%%%%%%%%%%%%%%%%%%%%%%%%%%%%%%%%%%%%%%%%%

\noindent \textit{Proof of Proposition \ref{prop:correctors}.} The statement \eqref{conv:xhi} is quite classical. The reader is referred to \cite[Ch. 2]{olla} for an insight of the method and to \cite[Prop. 4.3]{rhodes:06} for a proof in a more general context. \qed

\noindent \textit{Proof of Proposition \ref{prop:correctors2}.}   In what follows, for each $i=1,\cdots,d$, $({\bfm \varphi}_n^i)_n$ stands for a sequence in $\mathcal{C}$ such that $D{\bfm \varphi}_n^i \to {\bfm \zeta}^i$ in $L^2(\Omega)^d$ as $n\to +\infty$.

Let us first focus on \eqref{varform}. Fix $X\in\R^d$ whose entries are denoted by $(X_i)_{1\leq i \leq d}$.  We have: $D(X_i{\bfm \varphi}_n^i)=X_iD{\bfm \varphi}_n^i\to X_i{\bfm \zeta}_i={\bfm \zeta}X$  in $L^2(\Omega)^d$ as $n\to +\infty$ and
\begin{align*}
X^*\bar{A}X=&\M\big[(X+{\bfm \zeta}X)^*{\bfm a}(X+{\bfm \zeta}X)\big]
=\lim_{n\to +\infty}\M\big[(X+D(X_i{\bfm \varphi}_n^i))^*{\bfm a}(X+D(X_i{\bfm \varphi}_n^i))\big]\\
\geq &\inf_{{\bfm \varphi}\in \mathcal{C}}\M[(X+D{\bfm
\varphi})^*{\bfm a}(X+D{\bfm\varphi})].
\end{align*}
Conversely, from Lemma \ref{orthol} below, we have: 
\begin{equation}\label{orthogo}
\forall Y\in \R^d,\quad \M[(Y+{\bfm
\zeta}Y)^*{\bfm a}{\bfm \zeta}X]=\lim_{n\to +\infty}\M[(Y+{\bfm
\zeta}Y)^*{\bfm a}D(X_i{\bfm \varphi}_n^i)]=0.
\end{equation}
The above relation and Lemma \ref{orthol} again yield $\,\M[(X+{\bfm \zeta}X)^*{\bfm a}(D{\bfm\varphi}-{\bfm \zeta}X)]=0$ for any ${\bfm \varphi}\in\mathcal{C}$. So, for every ${\bfm \varphi}\in\mathcal{C}$, we have:
\begin{align*}
\M\big[(X+D{\bfm \varphi})^*{\bfm a}(X+D{\bfm \varphi})\big] =&\M\big[(X+{\bfm \zeta}X+D{\bfm \varphi}-{\bfm \zeta}X)^*{\bfm a}(X+{\bfm \zeta}X+D{\bfm \varphi}-{\bfm \zeta}X)\big] \\
=&\M\big[(X+{\bfm \zeta}X)^*{\bfm a}(X+{\bfm \zeta}X)\big] +2\M\big[(X+{\bfm \zeta}X)^*{\bfm a}(D{\bfm \varphi}-{\bfm \zeta}X)\big] \\&+\M\big[(D{\bfm \varphi}-{\bfm \zeta}X)^*{\bfm a}(D{\bfm \varphi}-{\bfm \zeta}X)\big] \\
\geq &\M\big[(X+{\bfm \zeta}X)^*{\bfm a}(X+{\bfm \zeta}X)\big] 
\end{align*}
so that \eqref{varform} follows. By the way, \eqref{orthogo} proves that $\bar{A}$ also matches $\M[(\mathrm{I}+{\bfm \zeta}^*){\bfm a}]$.

Now we prove $\Lambda {\rm I}\leq \bar{A}$. Fix $X\in \R^d$. From \eqref{unif_ellip} and Cauchy-Schwarz's inequality, we get
\begin{align*}
X^*\bar{A}X&=\M\big[(X+{\bfm \zeta}X)^*{\bfm a}(X+{\bfm \zeta}X)\big]\geq\Lambda \M\big[|X+{\bfm \zeta}X|^2\big]\geq\Lambda \big|\M\big[X+{\bfm \zeta}X\big]\big|^2=\Lambda|X|^2,
\end{align*}
since $\M[{\bfm \zeta}X]=0$. The estimate $\Lambda{\rm I}\leq \bar{A}$ follows.

Now we prove that $\bar{\Gamma}=\M[({\rm I}+{\bfm \zeta}^*){\bfm \gamma}] $ coincides with the orthogonal projection
$\M_1[({\rm I}+{\bfm \zeta}^*){\bfm \gamma}] $. Remind that ${\bfm \gamma}$  can be rewritten as ${\bfm \gamma}={\bfm a}e_1%+\bar{{\bfm H}}e_1
$. %(see Section \ref{notations} for the definition of $\bar{{\bfm H}}$ and \eqref{defcoef}). 
So we just have to establish the relation
\begin{equation}\label{anci}%i) \,\,\,\M_1[({\rm I}+{\bfm \zeta}^*)\bar{{\bfm H}}e_1] =0,\quad ii)\,\,\,
 \M_1[({\rm I}+{\bfm \zeta}^*){\bfm a}e_1] =\M[({\rm I}+{\bfm \zeta}^*){\bfm a}e_1] .
 \end{equation}
\textit{Proof of \eqref{anci}}. %The i-th entry $ M_1[(e_i+{\bfm \zeta}e_i)^*\bar{{\bfm H}}e_1] $ of $\M_1[({\rm I}+{\bfm \zeta}^*)\bar{{\bfm H}}e_1] $ satisfies ($\delta$ denotes the Kronecker symbol)
%\begin{equation*}
%M_1[(e_i+{\bfm \zeta}e_i)^*\bar{{\bfm H}}e_1]  =\M_1[({\rm I}+{\bfm \zeta}^*)_{ij}D_k{\bfm
%H}_{jk}]=\lim_{n\to
%\infty}\M_1[(\delta_{ij}+D_j{\bfm \varphi}_n^i)D_k{\bfm
%H}_{jk}]
%\end{equation*}
%Since ${\bfm
%H}_{ij}=0$ if $i=1$ or $j=1$, the index $k$ varies from $2$ to $d$ in the latter expression. By using  Lemma \ref{simple} iii) below, we deduce
%$$\lim_{n\to
%\infty}\M_1[(\delta_{ij}+D_j{\bfm \varphi}_n^i)D_k{\bfm
%H}_{jk}]=-\lim_{n\to \infty}\M_1[D^2_{jk}{\bfm \varphi}_n^i{\bfm
%H}_{jk}]=0$$ because of the antisymmetry of ${\bfm H}$. i) follows.
%\noindent \textit{Proof of ii)}.  
Because of the ergodicity of the measure $\mu$ (2. of Definition \ref{medium}), we stress that a function ${\bfm \psi}\in L^2(\Omega,{\cal G}^*,\mu) $ invariant  under the
translations $\{\tau_x;x\in\R\times\{0\}^{d-1}\}$ must be constant and therefore satisfies $\M_1[{\bfm \psi}]=\M[{\bfm \psi}]$. So we just have to prove that the entries $\M_1[(e_i+{\bfm \zeta}e_i)^*{\bfm a}e_1]$  are invariant under the translations $\{\tau_x;x\in\R\times\{0\}^{d-1}\}$. To that purpose, we only need to check that $$\M\big[\M_1[(e_i+{\bfm \zeta}e_i)^*{\bfm a}e_1]D_1{\bfm
\varphi}\big]=0$$ for any $i=1,\dots,d$ and ${\bfm \varphi}\in{\cal C}$. By using Lemma \ref{simple} ii below, we get:
\begin{align*} 
\M\big[\M_1[(e_i+{\bfm \zeta}e_i)^*{\bfm a}e_1]D_1{\bfm
\varphi}\big]&=\M\big[(e_i+{\bfm \zeta}e_i)^*{\bfm a}e_1\M_1[D_1{\bfm
\varphi}]\big]=\M\big[(e_i+{\bfm \zeta}e_i)^*{\bfm a}e_1D_1\M_1[{\bfm
\varphi}]\big].
\end{align*}
Since $D_k\M_1[{\bfm \varphi}]=0 $ for $k=2,\dots,d$ (see Lemma \ref{simple} i), we have $e_1D_1\M_1[{\bfm \varphi}]=D\M_1[{\bfm \varphi}]$. We deduce
$$\M\big[\M_1[(e_i+{\bfm \zeta}e_i)^*{\bfm a}e_1]D_1{\bfm
\varphi}\big]=\M\big[(e_i+{\bfm
\zeta}e_i)^*{\bfm a}D\M_1[{\bfm \varphi}]\big].$$
Since $\M_1[{\bfm \varphi}]\in \mathcal{C}$ (Lemma \ref{simple} ii), the latter quantity is equal to 0 (Lemma \ref{orthol}) and we complete the proof. Note that the above computations also prove: $\bar{\Gamma}_{1}=\M[(e_1+{\bfm \zeta}e_1)^*{\bfm a}e_1]=\bar{A}_{11}\geq \Lambda$.\qed

\begin{lemma}\label{orthol}
The following relation holds:
\begin{equation}\label{orthox}
\forall X\in\R^d,\quad \forall{\bfm \psi}\in\H,\quad \M\big[(X+{\bfm \zeta}X)^*{\bfm a}D{\bfm
\psi}\big]=0.
\end{equation}
\end{lemma}

\noindent \textit{Proof.} Since ${\bfm b}_i=\frac{1}{2}D_k{\bfm a}_{ik}$, the weak form of the
resolvent equation \eqref{wfe} associated to ${\bfm f}={\bfm b}_i$ reads, for any ${\bfm \psi}\in\H$: $$
\lambda({\bfm u}_\lambda^i,{\bfm
  \psi})_2+(1/2)\big({\bfm a}_{jk}D_j{\bfm u}_\lambda^i,D_k{\bfm
  \psi}\big)_2=(1/2)(D_k{\bfm a}_{ik},{\bfm \psi})_2=-(1/2)\big({\bfm a}_{ik},D_k{\bfm
  \psi}\big)_2.$$ By letting $\lambda$ go to $0$ and by using \eqref{conv:xhi}, we obtain: $(1/2)\big({\bfm a}_{jk}{\bfm \zeta}_j^i,D_k{\bfm
  \psi}\big)_2=-(1/2)\big({\bfm a}_{ik},D_k{\bfm
  \psi}\big)_2$. We deduce $\M\big[(\delta_{ij}+{\bfm \zeta}^i_j){\bfm a}_{jk}D_k{\bfm
\psi}\big]=0$, which means nothing but
\begin{equation}\label{ortho}
\M[(e_i+{\bfm \zeta}e_i)^*{\bfm a}D{\bfm
\psi}]=0.
\end{equation} 
The result follows by linearity.\qed

\begin{lemma}\label{simple}
The projection operator $\M_1$ saisfies the following elementary properties:

i)  $\forall k=2,\dots,d$ and $\forall {\bfm \varphi}\in {\rm Dom}(D_k)$, \,\,$D_k\M_1[{\bfm \varphi}]=\M_1[D_k{\bfm \varphi}]=0 $, 

ii)  $\forall {\bfm \varphi}\in {\cal C} $, \,\,$\M_1[{\bfm \varphi}]\in{\cal C}$
and \,\,$\M_1[D_1{\bfm \varphi}]=D_1\M_1[{\bfm \varphi}]$,

iii) $\forall k=2,\dots,d$ and $\forall
{\bfm \varphi},{\bfm \psi}\in {\rm Dom}(D_k)$, \,\,$\M_1[D_k{\bfm \varphi}{\bfm
\psi}]=-\M_1[{\bfm \varphi}D_k{\bfm \psi}] $. 
\end{lemma}

\noindent \textit{Proof.} The properties i) and ii)  are easily derived from the identities $\M_1[T_x{\bfm \varphi}]=\M_1[{\bfm \varphi}] $ for
any $x\in \{0\}\times \R^{d-1}$, $T_x\M_1=\M_1T_x$ for any
$x\in\R\times\{0\}^{d-1}$, and $\M_1[{\bfm \psi}*\rho]=\M_1[{\bfm
\psi}]*\rho$ for any ${\bfm \psi}\in L^\infty(\Omega)$ and $\rho\in
C^\infty_c(\R^d)$. iii) results from i). Details are left to
the reader.\qed

%%%%%%%%%%%%%%%%%%%%%%%%%%%%%%%%%%%%%%%%%%%%%%%%%%%%%%%%%%%%%%%%%%%%%%%%%%%%%%%%%%%%%%%%%%%%%%%%%%%%%%%%%%%
\section{J-topology}\label{sec:jaku}

%%%%%%%%%%%%%%%%%%%%%%%%%%%%%%%%%%%%%%%%%%%%%%%%%%%%%%%%%%%%%%%%%%%%%%%%%%%%%%%%%%%%%%%%%%%%%%%%%%%%%%
We summarized below the main properties of the Jakubowski topology (J-topology) on the space $D([0,T];\R)$ (set of functions that are right-continuous with left-limits on $[0,T]$) and refer the reader to \cite{jakubowski} for further details and proofs. We denote by $\mathds{V}$ the set of functions $v:[0,T]\to \R$ with bounded variations. The J-topology is a sequential topology defined by
\begin{definition}
A sequence $(x_n)_n$ in $D([0,T];\R)$ converges to $x_0\in D([0,T];\R)$ if for every $\varepsilon>0$, one can find elements $(v_{n,\varepsilon})_{n\in\nat}\subset \mathds{V}$ such that 

1) for every $ n\in\nat$,\, $\sup_{[0,T]}|x_n-v_{n,\varepsilon}|\leq \varepsilon$,

2) $\forall f:[0,T]\to \R$ continuous, \,$\int_0^Tf(r)dv_{n,\varepsilon}(r)\to\int_0^Tf(r)dv_{0,\varepsilon}(r)$ as $n\to +\infty$.
\end{definition}

By gathering \cite[Th. 3.8]{jakubowski} and \cite[Th. 3.10]{jakubowski}, one can state:
\begin{theorem}\label{thjaku}
Let $(V_\alpha)_\alpha\subset D([0,T];\R)$ be a family of nondecreasing stochastic processes. Suppose that the family $(V_\alpha(T))_\alpha$ is tight. Then the family $(V_\alpha)_\alpha$ is tight for the J-topology. Moreover, there exists a sequence $(V_n)_n\subset (V_\alpha)_\alpha$, a nondecreasing right-continuous process $V_0$ and a countable subset $C\subset[0,T[$ such that for all finite sequence $(t_1,\dots,t_p)\subset[0,T]\setminus C$, the family $(V_n(t_1),\dots, V_n(t_p))_n$ converges in law towards $(V_0(t_1),\dots, V_0(t_p))_n$.
\end{theorem}

Equip the set $\mathds{V}^+_c([0,T];\R)$ of continuous nondecreasing functions on $[0,T]$ with the J-topology and $C([0,T];\R)$ with the sup-norm topology. We claim:
\begin{lemma}\label{vc}
Let $(V_n)_n$ be a sequence in $\mathds{V}^+_c$ converging for the J-topology towards $V_0\in \mathds{V}^+_c$. Then $(V_n)_n$ converges towards $V_0$ for the sup-norm topology. 
\end{lemma}

\noindent \textit{Proof.} This results from Corollary 2.9 in \cite{jakubowski} and the Dini theorem.\qed

\begin{lemma}\label{contmap}
The following mapping is continuous
$$(x,v)\in C([0,T];\R)\times \mathds{V}^+_c([0,T];\R)\mapsto \int_0^\cdot x_r\,dv(r)\in C([0,T];\R).$$
\end{lemma}

\noindent \textit{Proof.} This results from Lemma \ref{vc} and the continuity of the mapping 
$$(x,v)\in C([0,T];\R)\times \mathds{V}^+_c([0,T];\R)\mapsto \int_0^\cdot x_r\,dv(r)\in C([0,T];\R),$$ where both $C([0,T];\R)$ and $\mathds{V}^+_c([0,T];\R)$ are equipped with the sup-norm topology. The reader may find a proof of the continuity of the above mapping in the proof of Lemma 3.3 in \cite{ouknine} (remark that, in \cite{ouknine}, the S-topology  coincides on $C([0,T];\R) $ with the sup-norm topology).\qed

%%%%%%%%%%%%%%%%%%%%%%%%%%%%%%%%%%%%%%%%%%%%%%
\section{Proof of the tightness (Proposition \ref{prop:tightness})}\label{app:tight}
%%%%%%%%%%%%%%%%%%%%%%%%%%%%%%%%%%%%%%%%%%%%%%%

We now investigate the tightness of the process $X^\varepsilon$ (and $K^\varepsilon $). %To
%that purpose, it is sufficient to establish the tightness of the
%process
%\begin{equation}\label{bcomp}
%\frac{1}{\varepsilon}\int_0^t{\bfm b}(\tau_{X^\varepsilon_r/\varepsilon}\omega)\,dr +
%\int_0^t{\bfm \gamma}(\tau_{X^\varepsilon_r/\varepsilon}\omega)\,dK^\varepsilon_r.
%\end{equation}
Roughly speaking, our proof is inspired by \cite[Chap. 3]{olla} and is based on the Garsia-Rodemich-Rumsey inequality: 
\begin{proposition}\label{rumsey}{\bf (Garsia-Rodemich-Rumsey's
inequality).} Let $p$ and $\Psi$ be strictly increasing continuous
functions on $[0,+\infty[$ satisfying $p(0)=\Psi(0)=0$ and
$\lim_{t\to \infty}\Psi(t)=+\infty$. For given $T>0$ and $f\in
C([0,T];\R^d)$, suppose that there exists a finite $B$ such that;
\begin{equation}\label{GRRcond}
  \int_0^T\int_0^T\Psi\Big(\frac{|g(t)-g(s)|}{p(|t-s|)}\Big)\,ds\,dt\leq B<\infty.
\end{equation}
Then, for all $0\leq s \leq t \leq T$:\quad $  |g(t)-g(s)|\leq 8\int_0^{t-s}\Psi^{-1}(4B/u^2)\,dp(u)
$.
%\begin{equation}\label{GRR}
%  |g(t)-g(s)|\leq 8\int_0^{t-s}\Psi^{-1}(4B/u^2)\,dp(u).
%\end{equation}
\end{proposition}

To apply Proposition \ref{rumsey}, it is necessary to establish exponential bounds for the drift of $X^\epsilon$. Indeed, suppose that we can prove the following exponential bound: for every $ 0\leq s,t\leq T$
\begin{equation}\label{eqron2}
  \bar{\E}^{\varepsilon *}\Big[\exp\Big(\kappa\big|\int_s^t\big[\frac{1}{\varepsilon}{\bfm b}_j-\partial_{x_i}
V(X^\varepsilon_r){\bfm a}_{ij}\big](\tau_{X^\varepsilon_r/\varepsilon}\omega)]\,dr+\int_s^t{\bfm a}_{1j}(\tau_{X^\varepsilon_r/\varepsilon}\omega)\,dK^\varepsilon_r \big|\Big)\Big]\leq 2\exp\big(C\kappa^2(t-s)\big).
\end{equation} for some constant $C>0$ depending only on $\Lambda$ (defined in \eqref{unif_ellip}). Then we can apply  Proposition \ref{rumsey} as detailed in \cite[Ch. 3, Th 3.5]{olla} (set $ p(t)=\sqrt{t}$,
$\psi(t)=e^t-1$ and $\psi^{-1}(t)=\ln(t+1)$ in  Proposition \ref{rumsey})
to obtain
\begin{proposition}\label{tightrho}We have the following estimate of the modulus of
continuity
\begin{align}\label{prop64}
\bar{\E}^{\varepsilon *}\Big(\sup_{|t-s|\leq \delta;0\leq
s,t\leq
T}\big|\int_s^t\big[\frac{1}{\varepsilon}{\bfm b}_j-\partial_{x_i}
V(X^\varepsilon_r){\bfm a}_{ij}\big](\tau_{X^\varepsilon_r/\varepsilon}\omega)]\,dr+&\int_s^t{\bfm a}_{1j}(\tau_{X^\varepsilon_r/\varepsilon}\omega)\,dK^\varepsilon_r  \big|\Big)\leq
C\sqrt{\delta}\ln(\delta^{-1}),
\end{align}
for some constant $C$ that only depends on $T,\Lambda$.
\end{proposition}

We easily deduce the proof of Proposition \ref{prop:tightness}:  we first work under $ \bar{\P}^{\epsilon*}$.
Let us investigate the tightness of $X^\varepsilon$.  Observe that
$$X^{j,\varepsilon}_t\!=\!x_j+\int_0^t\big[\frac{1}{\varepsilon}{\bfm b}_j(\tau_{X^\varepsilon_r/\varepsilon}\omega)-\partial_{x_i}
V(X^\varepsilon_r){\bfm a}_{ij}(\tau_{X^\varepsilon_r/\varepsilon}\omega)\big]\,dr+\int_0^t{\bfm a}_{1j}(\tau_{X^\varepsilon_r/\varepsilon}\omega)\,dK^\varepsilon_r
%+\int_0^t({\bfm \gamma}_j-{\bfm a}_{1j})(\tau_{X^\varepsilon_r/\varepsilon}\omega)\,dK^\varepsilon_r
+\int_0^t{\bfm \sigma}_{ji}(\tau_{X^\varepsilon_r/\varepsilon}\omega)\,dB^{*i}_r.$$
The tightness of the martingale part follows from the boundedness of ${\bfm
\sigma}$ and the Kolmogorov criterion. The tightness of the remaining terms results from Proposition \ref{tightrho}. So $X^\epsilon$ is tight under $ \bar{\P}^{\epsilon*}$.

Let us now investigate the tightness of the family
$(K^\varepsilon)_\varepsilon$. From Lemma \ref{lem_invariant}, we have
$\bar{\E}^{\varepsilon *}[K^\varepsilon_T]=T\int_{\partial
D}e^{-2V(x)}\,dx$. Theorem \ref{thjaku} ensures that
$(K^\varepsilon)_\varepsilon$ is tight in $D([0,T];\R_+)$ (remind that $K^\varepsilon$ is increasing). 

To sum up, under $ \bar{\P}^{\epsilon*}$, the family $(X^\varepsilon,K^\varepsilon)_\varepsilon $
is tight in $C([0,T];\bar{D})\times D([0,T];\R_+) $ equipped with
the product topology. From \eqref{girsanov}, the family is tight in $C([0,T];\bar{D})\times D([0,T];\R_+) $ under $ \bar{\P}^{\epsilon}$. \qed

\vspace{2mm}
We have thus shown that the proof of Proposition \eqref{prop:tightness} boils down to establishing \eqref{eqron2}. So we now focus on the proof of \eqref{eqron2}. We want to adapt the arguments of \cite[Chap. 3]{olla}. However, the situation is more complicated due to the pushing of the local time when $X^\epsilon$ is located on the boundary $\partial D$. Our idea is to eliminate the boundary effects by considering first a truncated drift vanishing near the boundary: fix $\omega\in \Omega$ and a smooth function $\rho\in
C^\infty_b(\bar{D})$ satisfying
$\rho(x)=0$ whenever $x_1 \leq \theta$ for some $\theta>0$. For any
$\varepsilon>0$ and $j=1,\dots,d$, define the "truncated" drift
\begin{equation}\label{bj}
b_{\rho,j}^\varepsilon(x,\omega)=\frac{e^{2V(x)}}{2}\partial_{x_i}\big(e^{-2V(x)}{\bfm a}_{ij}(\tau_{x/\varepsilon}\omega)\rho(x)\big),
\end{equation}
which belongs to $C^\infty_b(\bar{D})$. Our strategy is the following: we derive  exponential bounds for the process $\int_0^tb_{\rho,j}^\varepsilon(X^\varepsilon_r,\omega)\,dr$. These estimates will depend on $\rho$. Then we shall prove that we can choose an appropriate sequence $(\rho_n)_n\subset C^\infty_b(\bar{D})$ preserving the exponential bounds and such that the sequence $\big(\int_0^tb_{\rho_n,j}^\varepsilon(X^\varepsilon_r,\omega)\,dr\big)_n$ converges as $n\to \infty$ towards the process involved in \eqref{eqron2}.

The exponential bounds are derived from a proper spectral gap of the operator ${\cal L}^{\varepsilon}_V\cdot+ \kappa b_{\rho,j}^\varepsilon \cdot$ with boundary condition ${\bfm \gamma}_i(\tau_{x/\varepsilon}\omega)\partial_{x_i}\cdot=0$ on $\partial D $. The particular truncation we choose in \eqref{bj} is fundamental to establish such a spectral gap because it preserves the "divergence structure" of the problem. Any other (and maybe more natural) truncation fails to have satisfactory spectral properties.

So we define the set $$C^{2,\varepsilon}_\gamma=\{f\in
C^2_b(\bar{D});{\bfm \gamma}_i(\tau_{x/\varepsilon}\omega)\partial_{x_i}f(x)=0\text{
for }x\in \partial D\}$$ and consider the Hilbert space $L^2(\bar{D};e^{-2V(x)}dx) $ equipped with its norm $|\cdot|_D$ and its inner product $(\cdot,\cdot)_D $. Given $\kappa>0 $ and $ \omega\in\Omega$, let $\psi_\omega^{\varepsilon,\kappa}\in C^\infty([0,T]\times \bar{D})\cap C_b^{1,2}$ be the unique solution of
$$\partial_t \psi_\omega^{\varepsilon,\kappa}={\cal L}^{\varepsilon}_V\psi_\omega^{\varepsilon,\kappa}+\kappa b_{\rho,j}^\varepsilon (\psi_\omega^{\varepsilon,\kappa}+1)\,\text{ on }[0,T]\times D,\quad{\bfm \gamma}_i(\tau_{\cdot/\varepsilon}\omega)\partial_{x_i}
\psi_\omega^{\varepsilon,\kappa}=0\,\,\text{ on }[0,T]\times\partial
D$$ with initial condition $\psi_\omega^{\varepsilon,\kappa}(0,\cdot)=0$ on $ \bar{D}$ (see Lemma \ref{ladyzh}). Then  $u_\omega^{\varepsilon,\kappa}=\psi_\omega^{\varepsilon,\kappa}+1\in C^{1,2}_b$ is
a bounded classical solution of the problem
\begin{equation}
\partial_t u_\omega^{\varepsilon,\kappa}={\cal L}^{\varepsilon}_Vu_\omega^{\varepsilon,\kappa}+
\kappa b_{\rho,j}^\varepsilon u_\omega^{\varepsilon,\kappa}\,\text{ on }[0,T]\times D,\quad
{\bfm \gamma}_i(\tau_{\cdot/\varepsilon}\omega\big)\partial_{x_i}
u_\omega^{\varepsilon,\kappa}=0\,\,\, \text{ on }[0,T]\times\partial
D,
\end{equation}
with initial condition $u_\omega^{\varepsilon,\kappa}(0,\cdot)=1$ on $\bar{D}$. Lemma \ref{FFK} and a straightforward calculation provide the
probabilistic representation
\begin{align*}
u_\omega^{\varepsilon,\kappa}(t,x)=&\E^{\varepsilon*}_x\Big[\int_0^t\kappa b^\varepsilon_{\rho,j}(X^\varepsilon_r,\omega)\exp\big(\int_0^r\kappa b_{\rho,j}^\varepsilon(X^\varepsilon_u,\omega)\,du\big)dr\Big]+1\\=&\E^{\varepsilon*}_x\Big[\exp\big(\kappa\int_0^tb_{\rho,j}^\varepsilon(X^\varepsilon_r,\omega)\,dr\big)\Big].
\end{align*}

\begin{lemma}\label{var:exp}
For each $\omega\in \Omega$, we have the estimate $
|u_\omega^{\varepsilon,\kappa}(t,\cdot)|^2_D\leq e^{2t
\pi_\omega^{\varepsilon,\kappa}}$ $(0\leq t \leq T)$, where
$\pi_\omega^{\varepsilon,\kappa}=\sup (\phi,{\cal
L}^{\varepsilon}_V\phi+\kappa b_{\rho,j}^\varepsilon\phi)_D$ and the $\sup$ is taken over $\{\phi\in
C^{2,\varepsilon}_\gamma,\,|\phi|_D^2=1\}$.
\end{lemma}

\noindent \textit{Proof.} We have:
\begin{align*}
\partial_t|u_\omega^{\varepsilon,\kappa}(t,\cdot)|^2_D  = & 2(u_\omega^{\varepsilon,\kappa},\partial_t
u_\omega^{\varepsilon,\kappa}(t,\cdot))_D\\  = &
2(u_\omega^{\varepsilon,\kappa},{\cal L}^{\varepsilon}_V
u_\omega^{\varepsilon,\kappa}+\kappa b_{\rho,j}^\varepsilon
u_\omega^{\varepsilon,\kappa}(t,\cdot))_D \leq
2\pi_\omega^{\varepsilon,\kappa}|u_\omega^{\varepsilon,\kappa}(t,\cdot)|^2_D.
\end{align*}
Since $|u_\omega^{\varepsilon,\kappa}(0,\cdot)|^2_D=1$, we complete the proof
with the Gronwall lemma.\qed

\begin{proposition}\label{propolla}
For any $\kappa>0$, $\varepsilon>0$ and $0\leq s,t\leq T$
$$\bar{\E}^{\varepsilon *}\big[\exp\big(\big|\kappa\int_s^tb_{\rho,j}^\varepsilon(X^\varepsilon_r,\omega)\,dr\big|\big)\big]\leq
2\exp\big(C\kappa^2(t-s)\big),$$ for some constant $C$ that only
depends on $\Lambda $ and $\sup_{x\in \bar{D}}|\rho(x)|$.
\end{proposition}

\noindent \textit{Proof.} By stationarity (resulting from Lemma \ref{lem_invariant}) and
Lemma \ref{var:exp}, we have
\begin{align}
\bar{\E}^{\varepsilon *} &
\big[\exp\big(\kappa\int_s^tb_{\rho,j}^\varepsilon(X^\varepsilon_r,\omega)\,dr\big)\big]
  \leq
\bar{\E}^{\varepsilon *}\big[\exp\big(\kappa\int_0^{t-s}b_{\rho,j}^\varepsilon(X^\varepsilon_r,\omega)\,dr\big)\big]\nonumber\\
=& \M_D^*[u_\omega^{\varepsilon,\kappa}(t-s,x)] \leq
\M|u_\omega^{\varepsilon,\kappa}(t-s,\cdot)|_D\leq
\M[\exp((t-s) \pi_\omega^{\varepsilon,\kappa})].\label{estpi}
\end{align}
It remains to estimate $\pi_\omega^{\varepsilon,\kappa}$. For any function $\phi\in
C^{2,\varepsilon}_\gamma$ such that $|\phi|_D^2=1$, we have
\begin{align*}
(b_{\rho,j}^\varepsilon(\cdot,\omega),\phi^2)_D &
=-({\bfm a}_{ij}(\tau_{\cdot/\varepsilon}\omega)\rho\phi,\partial_{x_i}\phi)_D \leq
\Lambda^{-1}\sup_{x\in \bar{D}}|\rho(x)||\partial_x \phi|_D=C|\partial_x
\phi|_D
\end{align*}
where we have set $C= \Lambda^{-1}\sup_{x\in \bar{D}}|\rho(x)|$. As
a consequence (the $\sup$ below are taken over  $\{\phi\in
C^{2,\varepsilon}_\gamma,\,|\phi|_D^2=1\}$)
\begin{align}
\pi_\omega^{\varepsilon,\kappa}& = \sup(\phi,{\cal
L}^{\varepsilon}_V\phi+\kappa
b_{\rho,j}^\varepsilon\phi)_D\nonumber\\
&\leq \sup\big\{-(1/2)(a_{ij}(\tau_{\cdot/\varepsilon}\omega)\partial_{x_i}\phi,\partial_{x_j}\phi)_D+\kappa
(b_{\rho,j}^\varepsilon(\cdot,\omega),\phi^2)_D\big\}\nonumber\\
& \leq \sup\big\{-(\Lambda/2)|\partial_x\phi|_{D}^2+\kappa
C|\partial_x\phi|_{D}\big\} \leq \kappa^2C^2/(2\Lambda).\label{pi}
\end{align}
The last inequality is obtained by optimizing the expression $-(\Lambda/2)x^2+\kappa
Cx$ with respect to the parameter $x\in\R$. Gathering \eqref{estpi} and \eqref{pi} then yields $$\bar{\E}^{\varepsilon *}\big[\exp\big(\kappa\int_s^tb_{\rho,j}^\varepsilon(X^\varepsilon_r,\omega)\,dr\big)\big]\leq\exp\big(C'\kappa^2(t-s)\big)$$
where $C'=\sup_{x\in \bar{D}}|\rho(x)|^2/(2\Lambda^3)$. We complete
the proof by repeating the argument for $-b_{\rho,j}^\varepsilon$
and using the inequality $\exp(|x|)\leq \exp(-x)+\exp(x)$.\qed

\vspace{2mm}
As explained above, we can replace $\rho $ in Proposition   \ref{propolla}   with an appropriate sequence $(\rho_n)_n\subset C^\infty_b(\bar{D})$ so as to make the sequence $\big(\int_0^tb_{\rho_n,j}^\varepsilon(X^\varepsilon_r,\omega)\,dr\big)_n$ converging as $n\to \infty$ towards the process involved in \eqref{eqron2}. Let us construct such a sequence. For each $n\in\nat^*$, let us consider the  piecewise affine function
$\rho_n:\bar{D}\to \R$ defined by: 
$$\rho_n(x)=0 \text{ if } x_1\leq \frac{1}{n}, \quad \rho_n(x)=n(x_1-\frac{1}{n})\,\,\text{ if }\,\,\frac{1}{n}\leq x_1\leq \frac{2}{n}, \quad \text{ and } 1 \text{ otherwise}.$$
Note that $\rho_n $ is  continuous and $\sup_{x\in \bar{D}}|\rho_n(x)|\leq 1$. With the help of a regularization procedure and Lemma \ref{lem_invariant}, one can prove that Proposition \ref{propolla} remains valid for $\rho_n$ instead of $\rho $, where 
\begin{equation}\label{bro}
\begin{split}
\int_0^tb_{\rho_n,j}^\varepsilon(X^\varepsilon_r,\omega)\,dr=&\int_0^t\big[\frac{1}{\varepsilon}{\bfm b}_j(\tau_{X^\varepsilon_r/\varepsilon}\omega)-\partial_{x_i}
V(X^\varepsilon_r){\bfm a}_{ij}(\tau_{X^\varepsilon_r/\varepsilon}\omega)\big]\rho_n(X^\varepsilon_r)\,dr
\\&+\int_0^t{\bfm a}_{ij}(\tau_{X^\varepsilon_r/\varepsilon}\omega)n\one_{[\frac{1}{n};\frac{2}{n}]}(X^\varepsilon_r)\,dr.
\end{split}
\end{equation}
You can obtain the latter expression by expanding \eqref{bj} with respect to the operator $ \partial_{x_i}$. 

Since $\sup_{x\in\bar{D}}|\rho_n(x)|= 1$ for each $n$, we deduce
\begin{equation}\label{eqron}
\forall n \in \nat,\forall 0\leq s,t\leq T,\quad \bar{\E}^{\varepsilon *}\big[\exp\big(\big|\kappa\int_s^tb_{\rho_n,j}^\varepsilon(X^\varepsilon_r,\omega)\,dr\big|\big)\big]\leq
2\exp\big(C\kappa^2(t-s)\big)
\end{equation}
 for some constant $C$ only depending on $\Lambda$. Now it remains to pass to the limit as $n\to \infty$ in \eqref{eqron}. From \cite[Prop 1.19]{cattiaux},
$\int_0^t{\bfm a}_{ij}(\tau_{X^\varepsilon_r/\varepsilon}\omega)n\one_{[\frac{1}{n};\frac{2}{n}]}(X^\varepsilon_r)\,dr$
converges a.s. towards
$\int_0^t{\bfm a}_{1j}(\tau_{X^\varepsilon_r/\varepsilon}\omega)\,dK^\varepsilon_r $ as
$n\to\infty$. Fatou's Lemma (as $n\to +\infty$) in \eqref{eqron} then yields \eqref{eqron2}. So we complete the proof.\qed
%where we have set $A^{\varepsilon,j}_t=\int_0^t\big[\frac{1}{\varepsilon}{\bfm b}_j(\tau_{X^\varepsilon_r/\varepsilon}\omega)-\partial_{x_i}
%V(X^\varepsilon_r){\bfm a}_{ij}(\tau_{X^\varepsilon_r/\varepsilon}\omega)\big]\,dr$.
%From now on, we thus
%consider a subsequence, still indexed
%$(X^\varepsilon,K^\varepsilon)_\varepsilon $, that converges in
%$C([0,T];\R^d)\times D([0,T];\R) $ towards $(\bar{X},\bar{K})\in
%C([0,T];\bar{D})\times D([0,T];\R_+)$. Note that $\bar{K}$ is
%necessary a nondecreasing process as a S-limit of a sequence of
%nondecreasing processes. We also consider a countable subset ${\cal
%S}\subset [0,T)$ such that for all finite sets $Q\subset
%[0,T]\setminus {\cal S}$,
%$(X^\varepsilon,K^\varepsilon)\xrightarrow[\text{Dist.}(Q)]{}(\bar{X},\bar{K})$.

\subsection*{Acknowledgements}The author wishes to thank G. Barles, P.L. Lions, S.Olla  and T. Souganidis for interesting discussions that led to the final version of the mansucript, and V. Vargas who suggested the use of replication methods.
%%%%%%%%%%%%%%%%%%%%%%%%%%%%%%%%%%%%%%%%%%%%%%%%%%%%%%%%%%%%%%%%%%%%%%%%%%%%%%%%%%%%%%%%%%%%%%%%%%%%%

\end{document}